\documentclass[a4paper]{article}

\title{Embedding mapping class groups into \\ a finite product of trees}
\author{David Hume}
\date{\today}

\newcommand{\set}[1]{\left\{#1\right\}}
\newcommand{\setcon}[2]{\left\{#1\ \left|\ #2\right.\right\}}
\newcommand{\norm}[1]{\left\lVert#1\right\rVert}
\newcommand{\abs}[1]{\left\lvert#1\right\rvert}
\newcommand{\setg}[2]{[\![#1,#2]\!]}
\newcommand{\R}{\mathbb{R}}

\newcommand{\Z}{\mathbb{Z}}
\newcommand{\N}{\mathbb{N}}
\newcommand{\m}{\medskip \\}
\newcommand{\h}{\hspace{3mm}}
\newcommand{\tu}{\textup}

\newcommand{\mcT}{\mathcal{T}}
\newcommand{\geo}[1]{\underline{#1}}
\newcommand{\ogeo}[1]{\overline{\underline{#1}}}

\newcommand{\cT}[1]{\mathcal{T}(#1)}
\newcommand{\cC}[1]{\mathcal{C}(#1)}
\newcommand{\lv}[1]{\textup{lv}(#1)}

\newcommand{\bY}{{\bf Y}}
\newcommand{\bH}{{\bf H}}

\usepackage{amssymb}
\usepackage{amsmath, amssymb, amsthm}
\usepackage{enumitem}
\usepackage{mhsetup, mathtools}
\usepackage{tikz}
\usepackage{MnSymbol}
\usepackage{float}

\newtheorem{bthm}{Theorem}
\newtheorem{bcor}[bthm]{Corollary}
\newtheorem{thm}{Theorem}[section]
\newtheorem{cor}[thm]{Corollary}
\newtheorem{prop}[thm]{Proposition}

\newtheorem{lem}[thm]{Lemma}
\newtheorem*{prop*}{Proposition}
\newtheorem{defn}[thm]{Definition}
\newtheorem*{def*}{Definition}

\begin{document}
\maketitle

\begin{abstract}
\noindent We prove the equivalence between a relative bottleneck property and being quasi-isometric to a tree-graded space. As a consequence, we deduce that the quasi-trees of spaces defined axiomatically by Bestvina-Bromberg-Fujiwara are quasi isometric to tree-graded spaces. Using this we prove that mapping class groups quasi-isometrically embed into a finite product of simplicial trees. In particular, these groups have finite Assouad-Nagata dimension, direct embeddings exhibiting $\ell^p$ compression exponent $1$ for all $p\geq 1$ and they quasi-isometrically embed into $\ell^1(\N)$. We deduce similar consequences for relatively hyperbolic groups whose parabolic subgroups satisfy such conditions.

In obtaining these results we also demonstrate that curve complexes of compact surfaces and coned-off graphs of relatively hyperbolic groups admit quasi-isometric embeddings into finite products of trees.
\end{abstract}

\section{Introduction}

One of the most intensively studied classes of finitely generated groups are mapping class groups of compact surfaces due to their close connections with geometry, topology and group theory and their similarities with lattices in higher rank semisimple Lie groups and Out($F_n$). See for instance \cite{IvMCG,FarbMargMCG} and references therein.

In \cite{Gr87}, Gromov introduced relatively hyperbolic groups as a generalisation of hyperbolic groups. The class of relatively hyperbolic groups includes: hyperbolic groups, amalgamated products and HNN-extensions over finite subgroups, fully residually free (limit) groups \cite{Da03,Al05} - which are key objects in solving the Tarski conjecture \cite{Se01,KM10}, geometrically finite Kleinian groups and fundamental groups of non-geometric closed $3$-manifolds with at least one hyperbolic component \cite{Da03}. Mapping class groups are not relatively hyperbolic in general \cite{BDM09}.

Mapping class groups and relatively hyperbolic groups have been studied extensively from both algebraic and geometric perspectives. The goal of this paper is to consider these groups from the viewpoint of their quasi-isometric embeddings into finite products of (locally infinite) simplicial trees and coarse embeddings into $\ell^p$ spaces.
Many finitely generated groups are already known to admit quasi-isometric embeddings into a finite product of trees: hyperbolic, Coxeter, right-angled Artin and virtually special groups are all examples \cite{BDS07,DJ99,DJ00,HW08}. By contrast: the discrete Heisenberg group, Thompson's group and wreath products of infinite finitely generated groups admit no such embedding \cite{Pa01}. 

A natural metric generalisation of a tree is a quasi-tree, a geodesic metric space which is quasi-isometric to a tree. This important class of hyperbolic spaces is characterised by Manning's bottleneck property, \cite{Ma05}:

A geodesic metric space $X$ satisfies the \emph{bottleneck property} (BP) if and only if there is some constant $\Delta>0$ such that given any two distinct points $x,y\in X$ and some geodesic $\geo{g}$ from $x$ to $y$ with midpoint $m$, every path from $x$ to $y$ in $X$ intersects $B(m;\Delta)=\setcon{z\in X}{d_X(z,m)<\Delta}$.

When one considers relatively hyperbolic spaces (asymptotically tree-graded spaces in the sense of \cite{DS05}) the natural analogue of a tree is a tree-graded space. We recall that a geodesic metric space $X$ is tree-graded with respect to a collection of pieces $\setcon{X_i}{i\in I}$ if and only if each $X_i$ is closed and geodesic, $\abs{X_i\cap X_j}\leq 1$ whenever $i\neq j$ and any simple geodesic triangle is contained in a piece.

In this paper we define a \emph{relative bottleneck property} (cf. Definition \ref{RBP}) and prove an analogue of Manning's result.

\begin{bthm}\label{qTG} \h A geodesic metric space $X$ has the relative bottleneck property with respect to a collection of sets $\setcon{X_i}{i\in I}$ if and only if it is quasi-isometric to a space $\cT{X}$ which is tree-graded with pieces $\setcon{\mcT_j}{j\in J}$ where each $\mcT_j$ is either a point or is $(K,C)$ quasi-isometric to some $X_i$, where $K$ and $C$ are independent of $i$ or $j$.
\end{bthm}

Our key examples of spaces satisfying the relative bottleneck property (in a non--trivial way) are the quasi-trees of spaces defined by Bestvina-Bromberg-Fujiwara \cite{BBF10}. In this paper it is shown that mapping class groups quasi-isometrically embed into a finite product of quasi-trees of spaces, so, in particular they have finite asymptotic dimension. These spaces, denoted in this paper by $\cC{\bY}$, are constructed from a collection of spaces $\setcon{\cC{Y}}{Y\in\bY}$ - curve complexes in the mapping class group case, hence the choice of notation. The techniques in that paper have since been used to study embeddings of relatively hyperbolic groups into products of trees \cite{MS12}, where the collection of spaces consists of cosets of peripheral subgroups, which are not hyperbolic in general so lie outside the analysis conducted in \cite{BBF10}.

From Theorem \ref{qTG} we deduce that such a quasi-tree of spaces is quasi-isometric to a tree-graded space $\cT{\bY}$, with pieces which are either points or uniformly quasi-isometric to some $\cC{Y}$. From this we deduce several consequences for mapping class groups of compact surfaces and relatively hyperbolic groups. 

\begin{bcor}\label{MCGtrees} \h Mapping class groups of compact surfaces quasi-isometrically embed into a finite product of simplicial (but locally infinite) trees. In particular, they have finite Assouad-Nagata dimension, can be quasi-isometrically embedded into $\ell^1(\N)$ and, for each $p\in (1,\infty)$, admit explicit embeddings into $\ell^p$ spaces which exhibit compression exponent $1$.
\end{bcor}
The first two of these are consequences of the embedding into a product of trees but the third is more subtle and builds on the work in \cite{Hu11}. This corollary was previously only known in low complexity cases, where the mapping class group is virtually free, see for instance \cite{Be04}. A space with finite asymptotic dimension admits a coarse embedding into a Hilbert space, so mapping class groups satisfy the strong Novikov and coarse Baum-Connes conjectures \cite{BBF10,HR00, Yu00}. The Novikov conjecture had already been established independently by work of Hamenst\"{a}dt, Kida and Behrstock-Minsky. Kida, moreover, proves that mapping class groups are exact and hence have Yu's property (A) \cite{Ki08, Ha09, BM11}.

The $\ell^p$ compression exponent of a countable metric space $X$, $\alpha^*(X)$ - introduced in \cite{GK04} to quantify coarse embeddability - is the supremum over all $\alpha\in[0,1]$ such that there is some $C>0$ and a Lipschitz map $\phi:X\to\ell^p(\N)$ satisfying $\norm{\phi(x)-\phi(y)}_p\geq Cd(x,y)^\alpha$ for all $x,y\in X$. Compression exponents are closely linked to Yu's property (A) and amenability; and to the speed of random walks \cite{GK04,NP08}.

Assouad-Nagata dimension is a linearly-controlled version of Gromov's notion of asymptotic dimension \cite{As82,Gr93}. The Assouad-Nagata dimension of a space bounds the topological dimension of asymptotic cones \cite{DH08} and finite Assouad-Nagata dimension guarantees certain Lipschitz extension properties and $\ell^p$ compression exponent $1$ for all $p\geq 1$ \cite{BDHM, LS05, Gal}.
 
 We obtain similar consequences for relatively hyperbolic groups.

\begin{bcor}\label{RHdim} \h If $G$ is a finitely generated group, which is hyperbolic relative to a collection of subgroups $\setcon{H_i}{i\in I}$ then
\begin{itemize}
 \item $G$ has finite Assouad-Nagata dimension if and only if each $H_i$ does.
 \item $G$ can be quasi-isometrically embedded into $\ell^1(\N)$ if and only if each $H_i$ can,
 \item for each $p\geq 1$, the compression exponent $\alpha^*_p(G)=\min\setcon{\alpha^*_p(H_i)}{i\in I}$.
\end{itemize}
\end{bcor}
  The first of these was previously known for asymptotic dimension \cite{Os05}, the other two are generalisations of results contained in \cite{MS12,Hu11} respectively.

\subsection*{Plan of the paper}

Section \ref{introRBP} gives the precise definition of the relative bottleneck property and proves that it is satisfied by all quasi-trees of spaces constructed from the axiomatisation in \cite{BBF10}. We also prove that the property is a quasi-isometry invariant, which completes the reverse implication of Theorem \ref{qTG}. Section \ref{construction} gives the construction of a tree-graded space $\cT{X}$ from a space $X$ satisfying the relative bottleneck property and in section \ref{pfthm1} we prove that $\cT{X}$ is quasi-isometric to $X$ completing the forwards implication of Theorem \ref{qTG}. The final section (\ref{Conseq}) gives the full proof of Corollaries \ref{MCGtrees} and \ref{RHdim}.

\subsection*{Acknowledgements}

The author wishes to thank Cornelia Dru\c{t}u for many helpful conversations on the material in this paper, John Mackay and Alessandro Sisto for pointing out the more general version of their Theorem and Nicholas Loughlin for the introduction to tikZ. The author is grateful for the support of the EPSRC through a D.\! Phil.\! student grant and the grant ``Geometric and analytic aspects of infinite groups''. The paper was written while the author was a D.\! Phil.\! student at the University of Oxford. The author is indebted to the referee(s) who made many suggestions which greatly improved the clarity of the paper.

\section{Relative Bottleneck Property} \label{introRBP}
  In this section we introduce the relative bottleneck property, and state a key ``thickening'' lemma which will be used repeatedly throughout the paper. Also in this section we prove that the relative bottleneck property is a quasi-isometry invariant, and give the two key examples of spaces satisfying this property, tree-graded spaces and quasi-trees of spaces as defined in \cite{BBF10}.
  
We begin with the definition and terminology we will use during the paper.

\begin{defn}\label{RBP} \h {\bf Relative Bottleneck Property} \tu{(}cf. \tu{\cite{Ma05}}\tu{)}

Let $(X,d_X)$ be a geodesic metric space, let $\mathcal{X}=\setcon{X_i}{i\in I}$ be a collection of subsets of $X$ with $\bigcup_i X_i=X$ and let $M>0$. 

We say that the triple $(X,\mathcal{X},M)$ satisfies the \textbf{relative bottleneck property} \tu{(RBP)} if for each $i,j\in I$ with $i\neq j$ there is a tuple $I_{i,j}=(i=i_0,i_1,\dots,i_s=j)$ and for all $r\in\set{0,\dots s-1}$ there is some point $w_r\in X_{i_r}\cap X_{i_{r+1}}$ such that every path in $X$ from $X_i$ to $X_j$ intersects each of the open balls $B(w_r;M)\coloneqq\setcon{x\in X}{d_X(w_r,x)<M}$.
\end{defn}
  The following figure presents an idealised view of this definition. The focus of section \ref{sect:thickening} is to justify the extent to which this is a valid approximation.
\begin{figure}[H]
 \centering
 \begin{tikzpicture}[xscale=1, yscale=1, vertex/.style={draw,fill,circle,inner sep=0.3mm}]
\clip (-4.5,-2.5) rectangle (4.5,2.5);
{
\draw[black!50!white, very thin]
  (-4,1.65) circle (1cm);
 \draw[black!50!white, very thin]
  (-4.6,-1.1) circle (1cm);
 \draw[black!50!white, very thin]
  (5,-1.1) circle (1cm);
 \draw[black!50!white, very thin]
  (3.2,-2.4) circle (1cm);
 \draw[black!50!white, very thin]
  (-1.9,-1.6) circle (1cm);
 \draw[black!50!white, very thin]
  (0.4,-2.4) circle (1cm);
 \draw[black!50!white, very thin]
  (-1.1,1.7) circle (0.3cm);
 \draw[black!50!white, very thin]
  (-0.8,2) circle (0.3cm);
 \draw[black!50!white, very thin]
  (-1.35,2.1) circle (0.3cm);
 \draw[black!50!white, very thin]
  (-1.5,2.8) circle (0.5cm);
 \draw[black!50!white, very thin]
  (0.65,2.8) circle (1.4cm);
 \draw[black!50!white, very thin]
  (5.8,4.55) circle (4cm);

\filldraw[draw=black, fill=black!50!white, fill opacity=0.2, very thin]
 (-3,0) circle (1 cm);
 \path (-3,0) node[] {$X_i$};
 \filldraw[draw=black, fill=black!50!white, fill opacity=0.2, very thin] 
  (-1.2,0.5) circle (1 cm);
 \path (-1.2,0.5) node[] {$X_{i_1}$};
  \filldraw[draw=black, fill=black!50!white, fill opacity=0.2, very thin] 
  (0.4,-0.5) circle (1 cm);
 \path (0.4,-0.5) node[] {$X_{i_2}$};
   \filldraw[draw=black, fill=black!50!white, fill opacity=0.2, very thin] 
  (2.8,0) circle (1.5 cm);
 \path (2.8,0) node[] {$X_j$};

 \filldraw[fill=black, fill opacity=0.7]
 (-2.1,0.25) circle (0.4 cm);
 \node[vertex, color=white]
(c) at ( -2.1, 0.25) {};
\path (c) node[below, color=white] {$w_0$};
 \filldraw[fill=black, fill opacity=0.7]
 (-0.4,0) circle (0.4 cm);
 \node[vertex, color=white]
(d) at ( -0.4, 0) {};
\path (d) node[below, color=white] {$w_1$};
\filldraw[fill=black, fill opacity=0.7]
 (1.35,-0.25) circle (0.4 cm);
 \node[vertex, color=white]
(e) at ( 1.35, -0.25) {};
\path (e) node[below, color=white] {$w_2$};
}
\end{tikzpicture}

 \caption{The relative bottleneck property}\label{figRBP}
\end{figure}
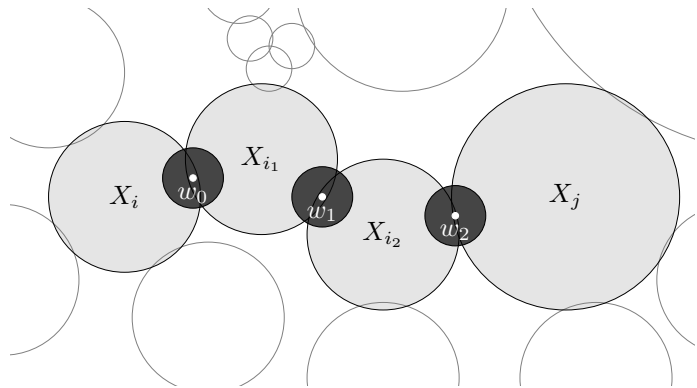

  Given a triple $(X,\mathcal{X},M)$ which satisfies (RBP) we call the elements of $\mathcal{X}$ \textbf{pieces}, and we call $M$ the \textbf{relative bottleneck constant}. In the remainder of the paper, when we state ``$(X,\mathcal{X},M)$ satisfies (RBP)'' we assume that for each pair of distinct elements $i,j\in I$ the tuple $I_{i,j}$ has been fixed and that a choice of suitable elements $w_r$ has been made and fixed. Our definition does not assume that $I_{j,i}$ is the tuple $I_{i,j}$ in reverse order, indeed this is emphatically not the case in Proposition \ref{exTG}.
  
Given two distinct pieces $X_i,X_j$ we define $W_{i,j}\coloneqq \setcon{w_r}{r=0,\dots,s-1}$ to be the set of \textbf{bottleneck points} from $X_i$ to $X_j$ and call the balls $B(w_r;M)$ the \textbf{bottlenecks} from $X_i$ to $X_j$.
  
  \medskip
  As a sample of the techniques used in this paper we now present a simple consequence of the above definition.

\begin{lem}\label{lem:distantpaths} Suppose $(X,\mathcal{X},M)$ satisfies (RBP) and let $X_i,X_j\in\mathcal{X}$. If there exist two paths $P,P'$ which start in $X_i$ and finish in $X_j$ with $d_X(P,P')\geq 2M$, then $i=j$.
\end{lem}
\begin{proof} Suppose for a contradiction that $i\neq j$. By definition there exists some bottleneck $B=B(w;M)$ from $X_i$ to $X_j$ with $B\cap P,B\cap P'\neq\emptyset$. Thus, $d_X(P,P')\leq d_X(P,w)+d_X(w,P')<2M$ which is a contradiction.
\end{proof}

\subsection{Thickening pieces}\label{sect:thickening}

In this section we present a construction which allows us to assume that the pieces in a space satisfying (RBP) are robustly path-connected in some sense. We will make this precise shortly. 

Notice that we do not even assume in the definition that the pieces $X_i$ are connected.

\begin{lem}\label{qconv} \h If $(X,\mathcal{X},M)$ satisfies the relative bottleneck property, then each $X_i\in\mathcal{X}$ is $4M$--quasi-convex, in the following sense:

If $x,y$ lie in $N_C(X_i)\coloneqq\setcon{y\in X}{d_X(y,X_i)\leq C}$ and $\geo{g}$ is a geodesic from $x$ to $y$, then $\geo{g}$ is contained in the $2M+2\max\set{M,C}$--neighbourhood of $X_i$.
\end{lem}
As a shorthand we denote the set of all geodesics from $x$ to $y$ by $\setg{x}{y}$.

\begin{proof} Set $M'\coloneqq \max\set{M,C}$. Let $x',y'$ be the end points of any component of $\geo{g}$ outside $N_{M'}(X_i)$, so $d(x',X_i),d(y',X_i)=M'$ and let $m$ be the mid-point of this component. As pieces cover $X$, $m\in X_k$ for some $k\in I$. 

Fix some $\varepsilon>0$ and let $x'',y''\in X_i$ be points at distance at most $M'+\varepsilon$ from $x',y'$ respectively and choose $\geo{g_x},\geo{g_y}\in \setg{x''}{x'},\setg{y''}{y'}$ respectively. Notice that for all $t\in [0,M']$, $d_X(\geo{g_x}(t),X_i),d_X(\geo{g_y}(t),X_i)\geq t-\varepsilon$.

Now consider the following two paths $P_x,P_y$ from $m\in X_k$ to $X_i$: $P_x$ (resp. $P_y$) is obtained by following $\geo{g}$ from $m$ to $x'$ (resp. $y'$) and then following $\geo{g_x}$ to $x''$ (resp. $\geo{g_y}$ to $y''$).

Since $i\neq k$ there is some bottleneck point $w\in X_i$ so $P_x\cap B(w;M), P_y\cap B(w;M)\neq\emptyset$. Let $z_x\in P_x\cap B(w;M)$ and $z_y\in P_y\cap B(w;M)$.

By the above we see that $d_X(x'', y'')\leq d_X(x'',z_x)+d_X(z_x,w)+d_X(w,z_y)+d_X(z_y,y'') < 4M + 2\varepsilon$. As this can be done for all $\varepsilon>0$ we deduce that $\geo{g}$ is contained in the $2M+2M'$--neighbourhood of $X_i$.
\end{proof}

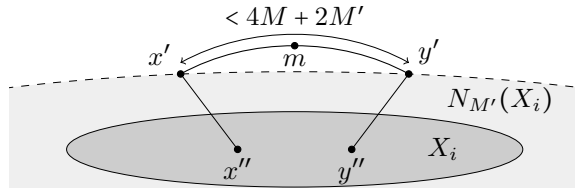
\begin{figure}[H]
 \centering
 \begin{tikzpicture}[xscale=1.5, yscale=0.25, vertex/.style={draw,fill,circle,inner sep=0.3mm}]
\clip (-2.5,-2) rectangle (2.5, 7.6);
{
\filldraw[draw=black, fill=black!20!white, fill opacity=0.3, very thin, dashed]
(0,0) circle (4.12 cm);
\filldraw[draw=black, fill=black!50!white, fill opacity=0.3, very thin]
 (0,0) circle (2cm);
 \draw[black, very thin]
 (-0.5,0) -- (-1,4);
 \draw[black, very thin]
 (0.5,0) -- (1,4);
 \draw[black]
 (-1,4) .. controls (-0.4, 6) and (0.4,6) .. (1,4);
 \draw[black, very thin, <->]
 (-1,4.6) .. controls (-0.4, 6.6) and (0.4,6.6) .. (1,4.6);

 \node[vertex]
(x) at ( -0.5, 0) {};
\path (x) node[below] {$x''$};
\node[vertex]
(y) at ( 0.5, 0) {};
\path (y) node[below] {$y''$};
 \node[vertex]
(x') at ( -1, 4) {};
\path (x') node[above left] {$x'$};
\node[vertex]
(y') at (1,4) {};
\path (y') node[above right] {$y'$};
\node[vertex]
(m) at (0,5.5) {};
\path (m) node[below] {$m$};
\path (1.3,0) node[] {$X_i$};
\path (1.8,3.8) node[below] {$N_{M'}(X_i)$};
\path (0,6.1) node[above] {$<4M+2M'$};
}
\end{tikzpicture}

 \caption{Quasi convexity of pieces}\label{figqconv}
\end{figure}

We would like to be in a situation where there are no sets of small diameter (compared to $M$) which disconnect a piece $X_i$. No such claim is made in the definition, but a simple ``thickening'' of the space achieves this. The robustness of the resulting connectivity of pieces is parametrised by a constant $b$ and - crucially - the bottleneck constant of the thickened space does not depend on $b$.

\begin{prop}\label{tech} \h Let $(X',\setcon{X_i'}{i\in I},M/9)$ satisfy \tu{(RBP)}. For every $b>0$ there is some  $(X^b,\setcon{X^b_i}{i\in I},M)$ which satisfies \tu{(RBP)}, and a $(2b+1)$-onto isometric embedding $\phi_b:X'\to X^b$ such that the restriction of $\phi^b$ to each $X'_i$ defines a $(2b+1)$-onto $(1,\frac{8M}{9}+1)$ quasi-isometric embedding $\phi^b_i:X'_i\to X^b_i$.
Moreover,
\begin{itemize}
\item there is a point $e$ (which will become the basepoint) contained in a unique piece $X^b_e$,
\item given any metric ball $B$ and any $X^b_i$ such that $B\cap X^b_i$ has diameter bounded by $2b$, $X^b_i\setminus B$ is path-connected.
\end{itemize}
\end{prop}
\begin{proof} Fix some $b>0$. Each piece $X'_i\in\mathcal{X}'$ is $\frac{4M}{9}$ quasi-convex by Lemma \ref{qconv}, so the $\frac{4M}{9}$--neighbourhoods of $X'_i$ (which we will label $X''_i$) are connected. Moreover, $(X',\setcon{X''_i}{i\in I},M)$ satisfies the relative bottleneck property.

We then achieve the first additional claim by defining a new point $e$ and attaching it to a unique piece $X''_e$ by a line of length $1$ (this line is added to $X''_e$). The resulting space under this construction so far is $(1,1)$ quasi-isometric to the original with uniformly $(1,\frac{8M}{9}+1)$ quasi-isometric pieces and has \tu{(RBP)} with constant $M$.

Now to achieve the second additional property we make the following construction.

We define $X^b_i=X''_i\times [0,2b+1]$ with the supremum product metric where the interval is given the standard Euclidean metric. Then we set
\[
X^b=\left.\bigsqcup_{i\in I} X^b_i\right/\sim \ \ \tu{where} \ \ (x,a)\sim(y,b) \ \tu{iff}\ a=b=0\ \tu{and}\ x=y.
\]
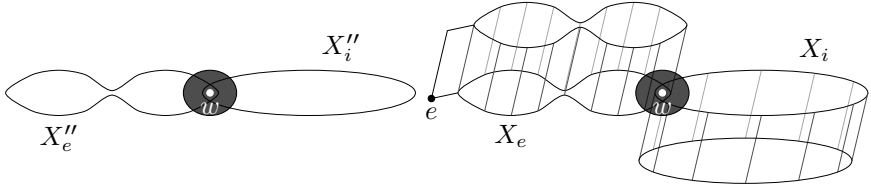
\begin{figure}[H]
 \centering
 \begin{tikzpicture}[xscale=0.7, yscale=0.15, vertex/.style={draw,fill,circle,inner sep=0.3mm}]
 \draw [black, xshift=-8.5cm, very thin] plot [smooth cycle] coordinates {(-2,0) (-1.5, 1.6) (-1,2) (-0.5, 1.6) (0,0.2) (0.5,1.6) (1,2) (1.5,1.6) (2,0) (1.5,-1.6) (1,-2) (0.5,-1.6) (0,-0.2) (-0.5,-1.6) (-1,-2) (-1.5,-1.6)};
 \draw [black, xshift=-8.5cm, very thin] (3.7,0) circle (2cm);

 \draw [black, very thin] plot [smooth cycle] coordinates {(-2,0) (-1.5, 1.6) (-1,2) (-0.5, 1.6) (0,0.2) (0.5,1.6) (1,2) (1.5,1.6) (2,0) (1.5,-1.6) (1,-2) (0.5,-1.6) (0,-0.2) (-0.5,-1.6) (-1,-2) (-1.5,-1.6)};
 \draw [black, very thin] (3.7,0) circle (2cm);
 
 \draw [black!70!white, opacity = 0.5, very thin]
 (-1.5,1.6)--(-1.2,7.6) (-1,2) -- (-0.7,8) (-0.5,1.6)--(-0.2,7.6) (0,0.2) -- (0.3,6.2) (0.5,1.6)--(0.8,7.6) (1,2) -- (1.3,8) (1.5,1.6)--(1.8,7.6);
 
 \foreach \r in {1,2,3,4,5}
 {
 \draw [black!70!white, opacity = 0.5, very thin]
 (3.7,0) ++(30*\r:2cm) -- ++(-0.3,-6);
 }

 \draw [black!70!white, very thin]
 (-2,0) -- (-1.7, 6) (-1.5,-1.6)--(-1.2,4.4) (-1,-2) -- (-0.7,4) (-0.5,-1.6)--(-0.2,4.4) (0,-0.2) -- (0.3,5.8) (0.5,-1.6)--(0.8,4.4) (1,-2) -- (1.3,4) (1.5,-1.6)--(1.8,4.4) (2,0)--(2.3,6);
 \draw [black, very thin]
 (-2,0)--(-2.5,-0.5)--(-2.2,5.5)--(-1.7,6);
 \foreach \r in {1,2,3,4,5,6,7}
 {
 \draw [black!70!white, very thin]
 (3.7,0) ++(150+30*\r:2cm) -- ++(-0.3,-6);
 }
 \path (-9.5,-2) node[below] {$X''_e$};
 \path (-4.2,2) node[above] {$X''_i$};
 
 \draw [black, yshift=6cm, xshift=0.3cm, very thin] plot [smooth cycle] coordinates {(-2,0) (-1.5, 1.6) (-1,2) (-0.5, 1.6) (0,0.2) (0.5,1.6) (1,2) (1.5,1.6) (2,0) (1.5,-1.6) (1,-2) (0.5,-1.6) (0,-0.2) (-0.5,-1.6) (-1,-2) (-1.5,-1.6)};
 \draw [black, yshift=-6cm, xshift=-0.3cm, very thin] (3.7,0) circle (2cm);
 
 \path (-1,-2) node[below] {$X_e$};
 \path (4.7,2) node[above] {$X_i$};
 
\filldraw[fill=black, fill opacity=0.7]
 (-6.66,0) circle[x radius=0.5cm, y radius=2cm];
 \node[vertex, color=white]
(f) at ( -6.66, 0) {};
\path (f) node[below, color=white] {$w$};
\filldraw[fill=black, fill opacity=0.7]
 (1.84,0) circle[x radius=0.5cm, y radius=2cm];
 \node[vertex, color=white]
(g) at ( 1.84, 0) {};
\path (g) node[below, color=white] {$w$};
\node[vertex]
(e) at ( -2.5, -0.5) {};
\path (e) node[below] {$e$};
\end{tikzpicture}
 \caption{The process in Proposition \ref{tech}}\label{figtech}
\end{figure}

It is clear that $X^b_i$ cannot be disconnected by a metric ball of diameter at most $2b$ with centre inside $X^b_i$. A ball centred outside $X^b_i$ which intersects this piece in a set of diameter at most $2b$ completely misses $X''_i\times\set{2b+1}$ so any two points $(x''_1,r_1)$ and $(x''_2,r_2)$ can be connected via $(x''_1,2b+1)$ and $(x''_2,2b+1)$ taking geodesics in the $[0,2b+1]$ direction and using the fact that $X''_i$ is connected. Also, as pieces only meet when the component of $[0,2b+1]$ is $0$ we have not changed the constant $M$. 

The natural injection $\phi^b$ of $X'$ into $X^b$ is a $(2b+1)$-onto isometric embedding and the restriction of $\phi^b$ to $X_i$ is a $(2b+1)$-onto $(1,\frac{8M}{9}+1)$ quasi-isometric embedding $\phi^b_i:X'_i\to X^b_i$.
\end{proof}

For completeness we note that $b=15M$ suffices for all arguments in this paper.

\subsection{Quasi-isometry invariance}

Theorem \ref{qTG} implies that \tu{(RBP)} is a quasi-isometry invariant, however, this is a straightforward consequence of the definition given by the following proposition.

\begin{prop}\label{RBPqi} \h Let $(X,d_X),(Y,d_Y)$ be geodesic metric spaces. If $q:X\to Y$ is a $(K,C)$-quasi-isometry and $(X,\setcon{X_i}{i\in I},M)$ satisfies (RBP) then there exists a constant $M'=M'(M,K,C)$ such that $(Y,\setcon{Y_i}{i\in I},M')$ satisfies (RBP) where $Y_i\coloneqq N_{C}(q(X_i))$.
\end{prop}
\begin{proof} It is clear that $\bigcup_{i\in I}Y_i = Y$ as $q$ is $C$-onto.

Let $i,j\in I$ with $i\neq j$ and let $w_k\in W_{i,j}$ be a bottleneck point from $X_i$ to $X_j$. We compute the distance between $q(w_k)\in Y_k\cap Y_{k+1}$ and some path $P$ from $Y_i$ to $Y_j$ in $Y$.

The pre-image under $q$ of $P$ defines a subset of $X$ whose $C$ neighbourhood contains a path from $N_{KC+C}(X_i)$ to $N_{KC+C}(X_j)$. Hence, $N_{KC+2C}(q^{-1}(P))\cap B(w_k;M)\neq\emptyset$. Applying $q$ we see that $d_Y(P,q(w_k))\leq K(KC+2C+M)+C$. 
\end{proof}

\subsection{Examples}
The two key examples of spaces satisfying \tu{(RBP)} are tree-graded spaces and quasi-trees of spaces satisfying the axioms of \cite{BBF10}.

Recall \cite{DS05} that a collection of subsets $\mathcal{X}$ of a geodesic metric space $X$ is called a \textbf{tree-grading} if the following three conditions hold: each $Y\in\mathcal{X}$ is geodesic; given any $Y,Z\in\mathcal{X}$ either $Y=Z$ or $Y\not\subseteq Z\not\subseteq Y$ and $Y\cap Z$ contains at most one point; any simple loop in $X$ is entirely contained in some $Y\in\mathcal{X}$. Notice that each $Y\in \mathcal{X}$ is a convex subset of $X$.

As an example, if we equip $\R^3$ with the metric 
\[
 d((x,y,z),(x',y',z'))= \left\{
  \begin{array}{lll} 
  \norm{(y-y',z-z')}_2
  &
  \tu{if}
  &
  x=x'
  \medskip
  \\
  \abs{x-x'}+\norm{(y,z)}_2 + \norm{(y',z')}_2
  &
  \tu{if}
  &
  x\neq x'
  \end{array}
  \right.
\]
we see that $\setcon{\set{x}\times\R^2}{x\in\R}$ is a tree-grading of $(\R^3,d)$, where the planes $\set{x}\times\R^2$ are equipped with the standard Euclidean norm.

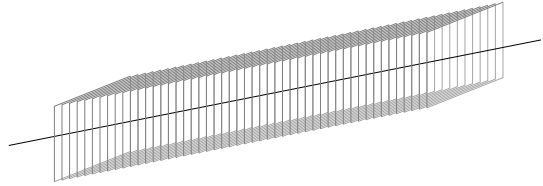
\begin{figure}[H]
 \centering
 \begin{tikzpicture}[xscale=1, yscale=1, vertex/.style={draw,fill,circle,inner sep=0.3mm}]
 \foreach \r in {1,...,50}{
 \draw[yslant=0.2,xshift=\r mm, black!50!white, opacity=0.7, very thin]
  (-0.5,-0.6) -- (0.5,-0.4) -- (0.5,0.6);
  }
 \draw[yslant=0.2,black, very thin] (-1,0)--(6,0);
 \foreach \r in {1,...,50}{
 \draw[yslant=0.2, xshift=\r mm, black!50!white, very thin]
  (0.5,0.6) -- (-0.5,0.4) -- (-0.5,-0.6);
  }
\end{tikzpicture}
 \caption{Accumulating pieces}\label{figtgR3}
\end{figure}
 Notice that distinct sets in this tree-grading are disjoint so we cannot hope to directly use the sets of a tree-grading to deduce (RBP). Instead we do the following.

\begin{prop}\label{exTG} \h Let $\setcon{X_i}{i\in I}$ be a tree-grading of a geodesic metric space $X$. Then $(X,\setcon{N_1(X_i)}{i\in I},2)$ satisfies \tu{(RBP)}.
\end{prop}

\begin{proof} 
Let $i,j\in I$, with $i\neq j$. We must choose a tuple $I_{i,j}=(i=i_0,\dots,i_s=j)$ and prove that there exist suitable bottleneck points.

Pick any geodesic $\geo{g}$ from some $x_i\in X_i$ to $x_j\in X_j\setminus X_i$. Choose $n_1\in\N$ maximal such that $\geo{g}(n_1)\in X_i$ and choose $i_1$ so that $\geo{g}(n_1+1)\in X_{i_1}$. If $d(x_i,x_j)<n_1+1$ then set $i_1=j$. Set $w_1=\geo{g}(n_1)\in X_i\cap N_1(X_{i_1})$.

Now choose $n_2\in\N$ maximal so that $\geo{g}(n_2)\in X_{i_1}$ and choose $i_2$ so that $\geo{g}(n_2+1)\in X_{i_2}$. If $d(x_i,x_j)<n_2+1$ then set $i_2=j$. Set $w_2=\geo{g}(n_2)\in X_{i_1}\cap N_1(X_{i_2})$.

Repeat this procedure until we define some $i_s=j$.

Now let $P$ be a path from $x'_i\in X_i$ to $x'_j\in X_j$ and suppose for a contradiction that $P\cap B(w_r;2)=\emptyset$ for some $r$.
Consider the loop $L'$ obtained by taking a geodesic from $x_i$ to $x'_i$, following $P$ to $x'_j$, taking a geodesic to $x_j$ and then $\geo{g}$ back to $x_i$.

$L'$ admits a simple subloop $L$ containing $w_r$ which is not a single point. In particular, this loop contains more than $1$ point in each of $X_{i_{r-1}}$ and $X_{i_{r}}$, so $L$ is not contained in a single piece. This contradicts the assumption that $\mathcal{X}$ is a tree-grading of $X$.
\end{proof}

Combined with Proposition \ref{RBPqi} this proves one direction of Theorem \ref{qTG}.

The second class of examples are quasi-trees of spaces defined axiomatically in \cite{BBF10}. 
The starting point of this construction is a collection of geodesic metric spaces $\setcon{\cC{Y}}{Y\in\bY}$ and for each ordered pair of distinct elements $Y,Z\in\bY$ a subset $\pi_Y(Z)\subseteq\cC{Y}$ of uniformly bounded diameter, which also satisfy a number of other axioms. 

We think of $\pi_Y(Z)$ as a `projection' of $\cC{Z}$ onto $\cC{Y}$. To give an example, if $G$ is a group which is hyperbolic relative to a subgroup $H$ we can define $\set{\cC{Y}}$ to be the collection of all cosets of $H$ in $G$ and consider $\pi_Y(Z)$ to be a closest point projection of the coset $Z$ onto the coset $Y$, that is, $\pi_Y(Z)=\setcon{g\in Y}{d(g,Z)=d(Y,Z)}$, where $d$ is some fixed word metric on $G$. 
The key example considered in \cite{BBF10} is projections between curve complexes of isotopic subsurfaces of a compact surface as defined by Mazur-Minsky \cite{MM00}.

These projections define functions $d^\pi_Y:\bY\setminus\set{Y}\times \bY\setminus\set{Y}\to\R$ given by $d^\pi_Y(X,Z)=\textrm{diam}(\pi_Y(X)\cup \pi_Y(Z))$, which are in some sense coarse pseudo--metrics. A technical point in the paper is the definition of functions $d_Y$ which differ from $d^\pi_Y$ by at most a uniform constant, but have more desirable properties.

Given spaces and projections satisfying suitable axioms, the first step is to build a \textbf{projection complex} $\mathcal{P}_K(\bY)$. This is a graph with vertex set $\bY$ where two vertices $X,Z$ span an edge if, for every $Y\not\in\set{X,Z}$, we have $d_Y(X,Z)\leq K$. Intuitively we imagine that in some hypothetical ambient space $\cC{X}$ and $\cC{Z}$ are `close'.

If $K$ is chosen to be sufficiently large, in comparison to other constants defined in the paper, then the projection complex $\mathcal{P}_K(\bY)$ is connected and is a quasi--tree. Given $X,Z\in\bY$ distinct, the set $\bY_K(X,Z)$ of all $Y\in\bY\setminus\set{X,Z}$ where $d_Y(X,Z)>K$ is precisely the set of internal vertices of a path from $X$ to $Z$ in $\mathcal{P}_K(\bY)$. We will assume from here on that a suitably large $K$ has been chosen so that all results from \cite{BBF10} may be applied.

Now we build an ambient space $\cC{\bY}_L$ from the disjoint union of the spaces $\setcon{\cC{Y}}{Y\in\bY}$ where for each edge $YZ\in\mathcal{P}_K(\bY)$ we connect each point in $\pi_Y(Z)$ to each point in $\pi_Z(Y)$ by a path of length $L$.
With this additional structure we can define the projection $\pi_Y(x)=\set{x}$ whenever $x\in\cC{Y}$. As a result one can define $d_Y(x,z)$ and $\bY_K(x,z)$ for all $x,z\in\cC{\bY}_L$ and all $Y\in\bY$. If $x\in\cC{X}$ and $z\in\cC{Z}$ then $\bY_K(x,z)\setminus \bY_K(X,Z)\subseteq \set{X,Z}$. If $\bY_K(x,z)=\emptyset$ then $d_{\cC{X}}(x,\pi_X(Z))\leq K$.

If $L$ is chosen appropriately compared to $K$ according to \cite[Lemma $4.2$]{BBF10} (and moreover $K\leq L\leq 2K$) then the subsets $\cC{Y}$ in $\cC{\bY}_L$ are totally geodesically embedded. We will assume from here on that such an $L$ has been chosen and replace $\cC{\bY}_L$ by $\cC{\bY}$.

The original inspiration for finding a relative bottleneck property is the following variation of \cite[Proposition $4.11$]{BBF10}.

\begin{prop}\label{prop:bbf4.11} Let $X, Z \in \bY$. If $Y\in\bY_K(X,Z)\cup\set{Z}$, then any path $P$ from $\cC{X}$ to $\cC{Z}$ in $\cC{\bY}$ contains a vertex $w$ such that
$d(w, \pi_Y(X)) < 7L$.
\end{prop}
\begin{proof} This follows immediately from \cite[Proposition $4.11$]{BBF10} whenever $P$ is a path from some $x\in\cC{X}$ to some $z\in\cC{Z}$ and $\bY_K(x,z)\neq\emptyset$.

Alternatively, if $\bY_K(x,z)=\emptyset$, then we see that $d_X(x,Z)\leq K$, so there is some point $z\in\pi_X(Z)$ such that $d(x,z)\leq K$. Now $XZ$ is an edge in $\mathcal{P}_K(\bY)$, so $d(x,\pi_Z(X))\leq K+L<2L$ and we are done.
\end{proof}

From this we prove that quasi-trees of spaces satisfy the relative bottleneck property, with respect to a collection of pieces uniformly quasi-isometric to the $\cC{Y}$.

\begin{prop}\label{exC(Y)} Let $\cC{\bY}$ be a quasi-tree of spaces satisfying the axioms of \tu{\cite{BBF10}}. Then $(X,\setcon{N_L(\cC{Y})}{Y\in\bY},10L)$ satisfies \tu{(RBP)}.
\end{prop}
\begin{proof} 
Given $X,Z\in\bY$ with $X\neq Z$ we define $I_{X,Z}$ to be the tuple $(X=Y_0,Y_1,\dots,Y_n=Z)$ where $\bY_K(X,Z)=\set{Y_1,\ldots,Y_{n-1}}$ and the $Y_i$ are arranged using the order property \cite[Theorem $3.3$(G)]{BBF10}.

To prove the proposition we will show that for each $i\in\set{1,\dots,n}$ there is a point $w_i\in \pi_{Y_i}(Y_{i-1})$ (and therefore in $N_L(\cC{Y_{i-1}})\cap N_L(\cC{Y_i})$) such that every path from $N_L(\cC{X})$ to $N_L(\cC{Z})$ contains a point in $B(w_i;10L)$.

By Proposition \ref{prop:bbf4.11} every path $P$ from $\cC{X}$ to $\cC{Z}$ contains some $x_i\in N_{8L}(\pi_{Y_i}(X))$ for all $i=1,\dots, m$. 

Next, we see that $diam(\pi_{Y_i}(X)\cup\pi_{Y_i}(Y_{i-1}))<L$ by the order and coarse equality properties \cite[Theorem $3.3$(G) and (B)]{BBF10}. Applying \cite[Lemma $4.2$]{BBF10} we see that for all $w_i$ in $\pi_{Y_i}(Y_{i-1})$ we have
\[
 d(P,w_i)\leq d(x_i,\pi_{Y_i}(X)) + diam(\pi_{Y_i}(X)\cup \pi_{Y_i}(Y_{i-1}))< 9L.
\]
Thus every path from $N_L(\cC{X})$ to $N_L(\cC{Z})$ contains a point in $B(w_i;10L)$.
\end{proof}

\subsection{Groups satisfying (RBP)}

  The relative bottleneck property is already well understood for finitely generated groups, via Stallings' Ends Theorem, which implies that Cay($G,S$) has (RBP) with respect to some subsets (in a non-trivial way) if and only if $G$ splits as an amalgam or HNN extension over a finite subgroup $G=A*_CB$ or $G=\tu{HNN}(A,C,\theta)$ (in a non-trivial way) \cite{St68,St71}. 
  
Any graph of groups decomposition induced in this way by the relative bottleneck property is accessible via results of Linnell \cite{Li83}, as the cardinality of subgroups over which we may amalgamate is uniformly bounded.

\section{Construction of the tree-graded space} \label{construction}
 From now on, we will assume - using Lemma \ref{tech} - that $X$ has \tu{(RBP)} with respect to a collection of pieces $\setcon{X_i}{i\in I}$ and a constant $M$ with a basepoint $e$ contained in a unique piece $X_e$ such that no metric ball which intersects $X_i$ in a set of diameter at most $2b$ disconnects $X_i$. As $M$ does not depend on $b$ from this point onwards we will assume that $b$ is sufficiently large, any $b\geq 15M$ will suffice. We fix such a $b$ for the remainder of the paper.

Our goal is to construct a suitable tree-graded space $\cT{X}$ which has the collection of pieces $\setcon{N_{4M}(X_i)}{i\in I}$.

For each $i\in I\setminus\set{e}$ we define $e_i\in X_i$ to be the point $w_0$ given by the bottleneck property such that all paths from $X_i$ to $X_e$ meet $B(e_i;M)$. Notice that $d(e,e_i)\leq d(e,X_i)+M$. We think of $e_i$ as a basepoint of $X_i$.

Our construction relies on organising pieces into strata parametrised by a (large) constant $R$. We fix a choice of $R\geq 160M$ for the remainder of the paper.

To this end we define a collection of strata $I^n\coloneqq\setcon{i\in I}{d(e,e_i)\leq nR}$ and set $I_n\coloneqq I^n\setminus I^{n-1}$. The level of $i$, $\lv{i}$ is the unique $n$ such that $i\in I_n$. By assumption $I^0=\set{e}$.

At this point we fix for each $X_i$ with $i\in I_{n+1}$ ($n\geq 0$) a geodesic $\geo{g_i}\in\setg{e_i}{e}$ and define $c_i$ to be the point on $\geo{g_i}$ at distance exactly $nR$ from $e$. We denote the reverse direction of a path $P$ by $\overline{P}$ and denote concatenation of paths by $P_1\circ P_2$, whenever the terminal point of $P_1$ agrees with the initial point of $P_2$.
\begin{figure}[H]
 \centering
 \begin{tikzpicture}[xscale=0.9, yscale=0.9, vertex/.style={draw,fill,circle,inner sep=0.3mm}]

\clip (-4.5,-2.5) rectangle (6.5,2.5);
{
\draw[black, very thin]
  (-4,1.65) circle (1cm);
 \draw[black!50!white, very thin]
  (-4.6,-1.1) circle (1cm);
 \draw[black!50!white, very thin]
  (5,-1.1) circle (1cm);
 \draw[black!50!white, very thin]
  (3.2,-2.4) circle (1cm);
 \draw[black!50!white, very thin]
  (-1.9,-1.6) circle (1cm);
 \draw[black!50!white, very thin]
  (0.4,-2.4) circle (1cm);
 \draw[black!50!white, very thin]
  (-1.1,1.7) circle (0.3cm);
 \draw[black!50!white, very thin]
  (-0.8,2) circle (0.3cm);
 \draw[black!50!white, very thin]
  (-1.35,2.1) circle (0.3cm);
 \draw[black!50!white, very thin]
  (-1.5,2.8) circle (0.5cm);
 \draw[black!50!white, very thin]
  (0.65,2.8) circle (1.4cm);
 \draw[black!50!white, very thin]
  (5.8,4.55) circle (4cm);
 \draw[xshift=7.5cm, yshift=-3.2cm, black!50!white, very thin]
  (-1.1,1.7) circle (0.3cm)  (-0.8,2) circle (0.3cm) (-1.35,2.1) circle (0.3cm);

\draw[black, very thin]
 (-3,0) circle (1 cm);
 \path (-3.3,0) node[] {$X_j$};
 \draw[black, very thin] 
  (-1.2,0.5) circle (1 cm);
 \path (-4,1.8) node[] {$X_{e}$};
  \draw[black, very thin] 
  (0.4,-0.5) circle (1 cm);
   \filldraw[draw=black, very thin, fill=black!30!white, fill opacity = 0.4] 
  (2.8,0) circle (1.5 cm);
 \path (2.8,0) node[] {$X_i$};
 
 \foreach \r in {1,...,7}
 {
 \draw[black!50!white, dashed]
 (-4,0.9) circle (1.8*\r cm);
 }
 }
 \node[vertex]
(g) at ( -4, 0.9) {};
\path (g) node[above left] {$e$};
\draw[black, thin] (-4,0.9)--(1.4,-0.25);
\node[vertex]
(j) at ( 1.4, -0.25) {};
\path (j) node[right] {$e_i$};
\node[vertex]
(k) at ( 1.28, -0.23) {};
\path (k) node[above left] {$c_i$};
\path (-0.8, 0.22) node[below] {$\geo{g_i}$};
 
 \end{tikzpicture}
 \caption{$R$-separated strata, in this example $i\in I_4$ and $j\in I_1$}\label{figstrata}
\end{figure}
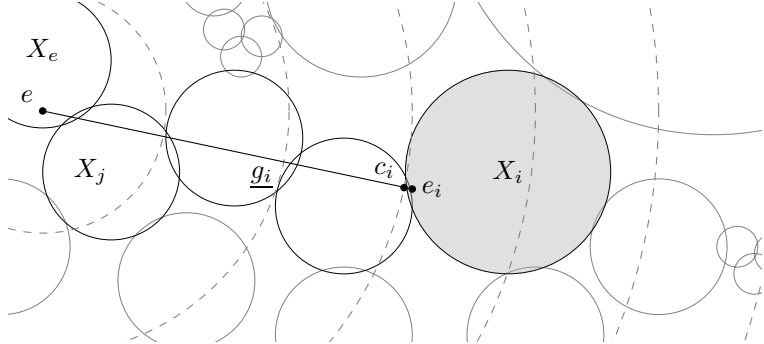
  The next two lemmas collect observations which will prove useful later. Before stating them we introduce some additional notation:
  
We call a path $P$ with endpoints $x$ and $y$ a $K$\textbf{-slack geodesic} if the length of $P$, $\abs{P}$ is bounded from above by $d_X(x,y)+K$.

\begin{lem}\label{10Mslack} For each $x\in N_{4M}(X_i)$ with $i\in I^n$ there is some $10M$-slack geodesic $\geo{q_x^i}$ from $x$ to $e$ which is contained in $N_{4M}(X_i)\cup B(e;nR)$.
\end{lem}
\begin{proof} If $x\in N_{4M}(X_i)$ with $i\in I^n$, then there is some $x'\in X_i$ with $d(x,x')\leq 4M$. We define the path $\geo{q_x^i}$ as the concatenation of some $\geo{g_1}\in\setg{x}{x'}$, $\geo{g_2}\in\setg{x'}{e_i}$ and $\geo{g_i}$. 

As $X_i$ is $4M$ quasi-convex by Lemma \ref{qconv} and $e_i\in B(e;nR)$, we see that $\geo{q_x^i}\subseteq N_{4M}(X_i)\cup B(e;nR)$. Every geodesic from $x'$ to $e$ passes within $M$ of $e_i$ by \tu{(RBP)}. Hence,
$\abs{q_x^i}\leq 4M + d(x',e) + 2M \leq d_X(x,e) + 10M$. \end{proof}
\begin{figure}[H]
 \centering
 \begin{tikzpicture}[xscale=1, yscale=1, vertex/.style={draw,fill,circle,inner sep=0.3mm}]
\clip (-4.5,-1.5) rectangle (4.5,2);
{
 \filldraw[draw=black!50!white, very thin, dashed, fill=black!10!white, fill opacity=0.3]
 (-4,0.9) circle (5.8cm);
 \filldraw[draw=black!50!white, fill=black!30!white, fill opacity=0.3, very thin] 
  (2.8,0) circle (1.5 cm);
 \path (2.8,0) node[] {$X_i$};

\filldraw[fill=black, fill opacity=0.6]
 (1.4,-0.25) circle (0.4 cm);
}
 \node[vertex]
(g) at ( -4, 0.9) {};
\path (g) node[above left] {$e$};

\draw[black, semithick] (3.2,1.8)--(3.2,1) -- (1.4,-0.25) -- (-4,0.9);
\node[vertex, color=white]
(j) at ( 1.4, -0.25) {};
\path (j) node[below, color=white] {$e_i$};
\node[vertex]
(k) at ( 3.2, 1.8) {};
\path (k) node[right] {$x$};
\node[vertex]
(l) at ( 3.2, 1) {};
\path (l) node[right] {$x'$};

\path (-0.9, 0.2) node[below] {$\geo{q_x^i}$};
\path (0.5, 1) node[above] {$B(e;nR)$};

\end{tikzpicture}
 \caption{Finding $10M$-slack geodesics}\label{figq_x^i}
\end{figure}
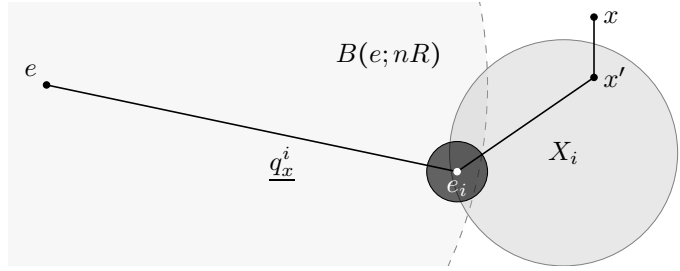

\begin{lem}\label{I^nlv} Let $i,j\in I$, with $i\neq j$. If $d_X(e_i,e)\geq d_X(e_j,e)$ then every path from $X_i$ to $X_j$ in $X$ contains a point in $B(e_i;4M)$.
\end{lem}
\begin{proof} Suppose there is a path $P$ from $x\in X_i$ to $y\in X_j$ which avoids the ball $B(e_i;M)$. If $d(e_i,e_j)\geq 2M$ then any geodesic in $\setg{e_j}{e}$ avoids this ball, and as $X_j$ has no small cut-sets there is a path from $y$ to $e$ also avoiding this ball (for instance extend such a path from $y$ to $e_j$ by $\geo{g_j}$) contradicting \tu{(RBP)}.

Now consider a path $P'$ of length at most $2M$ from $e_i$ to $e_j$, some point on this path lies in a bottleneck for paths between $X_i$ and $X_j$. By Lemma \ref{lem:distantpaths} any path $Q$ from $X_i$ to $X_j$ satisfies $d(Q,P')<2M$ and hence $Q\cap B(e_i;4M)\neq\emptyset$.
\end{proof}
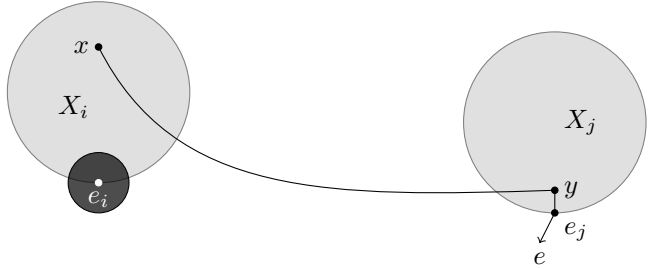
\begin{figure}[H]
 \centering
 \begin{tikzpicture}[xscale=1, yscale=1, vertex/.style={draw,fill,circle,inner sep=0.3mm}]
 \filldraw[draw=black!50!white, fill=black!30!white, fill opacity=0.4]
 (-3,0.4) circle (1.2 cm) (3,0) circle (1.2 cm);
 
\filldraw[fill=black, fill opacity=0.7]
 (-3,-0.8) circle (0.4 cm);
 \node[vertex, color=white]
(c) at ( -3, -0.8) {};
\path (c) node[below, color=white] {$e_i$};
\node[vertex]
(j) at ( 3, -1.2) {};
\path (j) node[below right] {$e_j$};
\draw[black] (-3,1) .. controls (-2,-1) and (0,-1) .. (3,-0.9)--(3,-1.2);
\draw[->, black] (3,-1.2)--(2.8,-1.6);
\path (2.8,-1.6) node[below] {$e$};
\node[vertex]
(k) at ( -3, 1) {};
\path (k) node[left] {$x$};
\node[vertex]
(l) at (3,-0.9) {};
\path (l) node[right] {$y$};
\path (-3,0.2) node[left] {$X_i$};
\path (3,0) node[right] {$X_j$};
\end{tikzpicture}

 \caption{Passing to lower levels when $d(e_i,e_j)\geq 2M$}\label{figlowlev}
\end{figure}
  One key element of this paper is deciding when pieces in the same level should have an immediate common ancestor in the tree-graded space. We introduce the following equivalence relation on each level $I_{n+1}$ to help determine this:

\begin{def*} 
Given $i,j\in I_{n+1}$ we write $i\sim j$ if and only if there exists some path $P$ from $X_i\setminus B(e;nR+11M)$ to $X_j\setminus B(e;nR+11M)$ with the property that $P\cap B(e;nR+11M)$ is either empty or is contained in $N_{4M}(X_k)$ for some $k\in I^n$.
\end{def*}
\begin{figure}[H]
 \centering
 \begin{tikzpicture}[xscale=1, yscale=1, vertex/.style={draw,fill,circle,inner sep=0.3mm}]
\clip (-5.5,-2.5) rectangle (5.5,1.5);
{
\filldraw[draw=black!30!white, dashed, very thin, fill=black!10!white, fill opacity =0.2]
 (0,-10) circle (9.2cm);
 \filldraw[draw=black!30!white, dashed, very thin, fill=black!10!white, fill opacity =0.2]
 (0,-10) circle (8.2cm);
\filldraw[draw=black!50!white, fill=black!30!white, fill opacity=0.4]
 (-2.5,0.1) circle (1.2 cm) (2.5,0) circle (1.2 cm) (0,-1.3) circle (1.2cm);

\filldraw[draw=black!50!white, fill=black!50!white, fill opacity=0.5]
 (-2.5,0.1) +(-45:1.2cm) circle (0.7 cm);
\node (A) at (-2.5,0.1) {};
\node[vertex]
(B) at ([shift={(-45:1.2 cm)}]A) {};
\path (B) node[below] {$e_i$};
\node (C) at (2.5,0) {};
\node[vertex]
(D) at ([shift={(-135:1.2 cm)}]C) {};
\path (D) node[below] {$e_j$};

\path (-1.7,-2.35) node[] {$B(e;nR)$};
\path (-2.6,-2.1) node[above] {$B(e;nR+11M)$};

\draw[black] plot [smooth] coordinates {(-2.5,0.3) (-1.5,-0.6) (-0.4,-1.3) (0.4,-1.3) (2,0)};

\path (-3.2,0.1) node[right] {$X_i$};
\path (3.2,0) node[left] {$X_j$};
\path (0,-2.5) node[above] {$X_{k}$};
}\end{tikzpicture}

 \caption{The equivalence $i\sim j$}\label{figisimj}
\end{figure}
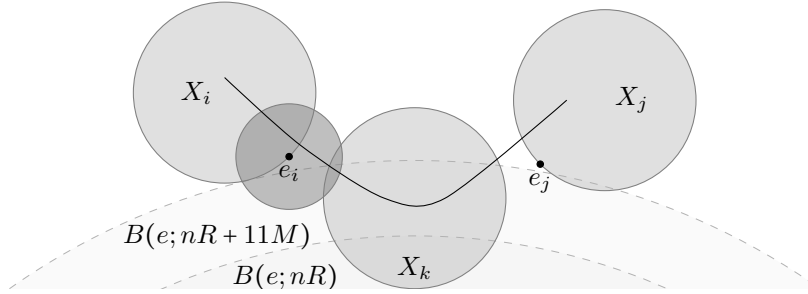
  Without loss, we may assume $d_X(e,e_i)\geq d_X(e,e_j)$, so such a path intersects the $4M$ ball around $e_i$ by Lemma \ref{I^nlv}. The fact that this does define an equivalence relation is not obvious so we provide a proof.

\begin{lem}\label{reln} \h $\sim$ is an equivalence relation.
\end{lem}
\begin{proof} Firstly we prove that $\sim$ is reflexive, for this it suffices to prove that whenever $i\in I_{n+1}$ then $X_i\setminus B(e;nR+11M)$ is not empty.

Any geodesic from a point $x$ in $X_i\cap B(e;nR+11M)$ to $e$ contains a point in $B(e_i;M)$. Since $d(e_i,e)\geq nR$ it follows that $d(x,e_i)\leq 13M$ so $X_i\cap B(e;nR+11M)$ has diameter at most $26M$. By Lemma \ref{tech}, and the fact that $b\geq 13M$, $X_i\setminus B(e;nR+11M)$ is path-connected, so $i\sim i$.

The fact that $\sim$ is symmetric is immediate. We now prove that it is transitive.

Suppose $i\sim j\sim l$ with $\abs{\set{i,j,l}}=3$. Since $i\sim j$ there is a path $P_1$ from $x_i\in X_i$ to $x_j^1\in X_j$ such that $P_1\cap B(e;nR+11M)$ is either empty or contained in some $N_{4M}(X_k)$ with $k\in I^n$. Likewise there is a path $P_2$ from $x_j^2\in X_j$ to $x_l\in X_l$ such that $P_2\cap B(e;nR+11M)$ is either empty or contained in some $N_{4M}(X_{k'})$ with $k'\in I^n$.

Using Lemma \ref{tech} as above we know that $x^1_j,x^2_j$ do not lie in $B(e;nR+11M)$ and hence we can find a path $Q\subset X_j$ between them which is disjoint from $B(e;nR+11M)$.

Hence the path $P_1\circ Q \circ P_2$ establishes that $i\sim l$ unless $k\neq k'$ and the sets $P_1\cap B(e;nR+11M)$, $P_2\cap B(e;nR+11M)$ are both non-empty. We now assume, for a contradiction, that this is the case. Consider the following two paths from $N_{4M}(X_k)$ to $N_{4M}(X_{k'})$.
\begin{itemize}
\item $\geo{g_k}\circ\ogeo{g_{k'}}$ (contained in $B(e;nR)$),
\item $P$ (avoids $B(e;nR+10M)$): follow $P_1$ from some point in  $N_{4M}(X_k)\setminus B(e;nR+10M)$ to $x_j^1$ then take $Q$ to $x^j_2$ and follow $P_2$ to a point in $N_{4M}(X_{k'})\setminus B(e;nR+10M)$.
\end{itemize}

 These paths are at distance at least $10M$ contradicting Lemma \ref{lem:distantpaths}. Hence, $k=k'$ which contradicts the initial assumption that $k\neq k'$.
  \end{proof}
\begin{figure}[H]
 \centering
 \begin{tikzpicture}[xscale=1, yscale=1, vertex/.style={draw,fill,circle,inner sep=0.3mm}]
\clip (-6,-2.7) rectangle (6,1.5);
{
\filldraw[draw=black!30!white, dashed, very thin, fill=black!10!white, fill opacity =0.2]
 (0,-10) circle (9.5cm);
 \filldraw[draw=black!30!white, dashed, very thin, fill=black!10!white, fill opacity =0.2]
 (0,-10) circle (8.5cm);
\filldraw[draw=black!50!white, fill=black!30!white, fill opacity=0.4]
 (-3.6,0.1) circle (1.2 cm) (0,0.3) circle (1.2 cm) (3.5,0) circle (1.2cm);

\filldraw[draw=black!50!white, fill=black!50!white, fill opacity=0.4]
 (-1.7,-1.3) circle (1 cm) (1.7,-1.4) circle (1 cm);
\node (A) at (-1.7,-1.3) {};
\node[vertex]
(B) at ([shift={(-60:1 cm)}]A) {};
\path (B) node[above left] {$e_k$};
\node (C) at (1.7,-1.4) {};
\node[vertex]
(D) at ([shift={(-120:1 cm)}]C) {};
\path (D) node[above right] {$e_{k'}$};
\draw[black, semithick, ->] (B) -- +(0.2,-0.4);
\draw[black, semithick, ->] (D) -- +(-0.2,-0.4);

\path (B) +(0.6,-0.35) node[right] {$\geo{g_k}\circ\ogeo{g_{k'}}$};

\path (-2.9,-2.5) node[] {$B(e;nR)$};
\path (4.3,-2.4) node[above] {$B(e;nR+11M)$};

\path (-1.7,-1.2) node[above] {$X_k$};
\path (1.8,-1.2) node[above] {$X_{k'}$};

\path (-3.6,0.5) node[] {$X_i$};
\path (0,0.8) node[] {$X_j$};
\path (3.5,0.5) node[] {$X_l$};

\draw[black, very thin] plot [smooth] coordinates {(-3.6,0.1) (-3,-0.6) (-2,-1.3) (-1.6,-1.3) (-1.2,-0.8) (-0.6,0.3)};
\node[vertex]
(E) at (-1.05,-0.55) {};
\path (-3.3,-0.2) node[right] {$P_1$};
\path (3.1,-0.2) node[left] {$P_2$};
\path (0.1,0.3) node[below] {$P$};
\draw[black, very thin] plot [smooth] coordinates {(0.8,0.3) (1.3,-0.8) (1.6,-1.3) (2,-1.3) (2.7,-0.65) (3.3,0)};
\node[vertex]
(F) at (1.18,-0.55) {};
\draw[black, semithick] (-1.05,-0.55)--(-0.6,0.3)--(0.8,0.3)--(1.18,-0.55);

}\end{tikzpicture}

 \caption{Transitivity of the relation $\sim$}\label{figsimtrans}
\end{figure}
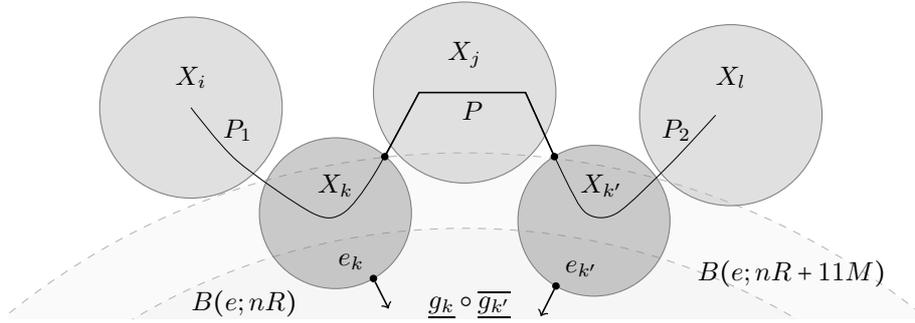
  In our construction of the tree-graded space $\cT{X}$ we will insist that whenever $i\sim j$, the pieces $N_{4M}(X_i)$ and $N_{4M}(X_j)$ are connected to the same piece in a lower level. The following two lemmas provide candidate pieces in lower levels.
  
\begin{lem}\label{highbottleneck} Let $i\in I_{n+1}$. There exists some $k\in I_{i,e}\cap I^n$ and some $w\in X_k\cap W_{i,e}$ such that $d_X(e,w)\geq nR-M$.
\end{lem}
\begin{proof} Let $k=i_m$ be the first coordinate of the tuple $I_{i,e}$ contained in $I^n$ and let $w$ be a bottleneck point contained in $X_{i_{m-1}}\cap X_{k}$. As $w\in X_{i_{m-1}}$, every geodesic from $w$ to $e$ meets $B(e_{i_{m-1}};M)$. Therefore, $d(w,e)\geq d(e_{i_{m-1}},e)-M\geq nR-M$, since $i_{m-1}\in I_{n+1}$.
\end{proof}
  
\begin{lem}\label{glue} \h Let $i\in I_{n+1}$. There exists some $k\in I^n$ such that 
\[
 \setcon{c_j}{j\in[i]}\subseteq N_{4M}(X_k).
\]
\end{lem}
  We will actually prove that this happens whenever $k\in I_{i,e}\cap I^n$ and there is some $w\in X_k\cap W_{i,e}$ with $d_X(e,w)\geq nR-M$. The existence of such a $k$ is given by Lemma \ref{highbottleneck}. It is not necessarily the case that $k\in I_{j,e}$ for all $j\in [i]$.

\begin{proof} Fix some $k$ with the above properties. We first prove that for every $j\in[i]$, there is some $w_j\in X_k$ such that $\geo{g_j}\cap B(w_j;M)\neq\emptyset$ and $d(w_j,e)\geq nR-M$. By assumption this holds when $i=j$.

Pick some $j\in[i]$ with $j\neq i$ and let $P$ be a path from $X_i$ to $X_j$ such that $P\cap B(e;nR+11M)\subseteq N_{4M}(X_{k'})$ for some $k'\in I^n$.

If $P\cap B(e;nR+11M)=\emptyset$ then either $\geo{g_j}\cap B\neq\emptyset$ or $P\cap B\neq\emptyset$, since every path from $X_i$ to $e$ meets $B$ and $B$ cannot disconnect $X_j$. In the first case we are done; in the second, we construct a path $P'$ from $X_j$ to $e$ such that all points $p\in P'$ satisfying $d(e,p)\in[nR,nR+10M]$ are contained in $X_k$. Recall that by Lemma \ref{tech}, $X_k\setminus B$ is path-connected whenever $B$ is a metric ball of radius $M$, so if $w'\in W_{j,e}$, then either $d(w',e)\leq nR+M$ or $d(w',e)\geq nR+9M$.

If $P\cap B(e;nR+11M)\neq\emptyset$ then, working as above, we obtain a path $P''$ from $X_j$ to $e$ via $e_k$ such that all points $p\in P''$ satisfying $d(e,p)\in[nR,nR+7M]$ are contained in $X_{k'}$. Therefore, for every $w'\in W_{j,e}$, either $d(w',e)\leq nR+M$ or $d(w',e)\geq nR+6M$.
\begin{figure}[H]
 \centering
 \begin{tikzpicture}[xscale=1, yscale=1, vertex/.style={draw,fill,circle,inner sep=0.3mm}]
\clip (-5.5,-3.5) rectangle (5.5,1.3);
{

\filldraw[draw=black!30!white, dashed, very thin, fill=black!10!white, fill opacity =0.2]
 (0,-10) circle (9.2cm);
 \filldraw[draw=black!30!white, dashed, very thin, fill=black!10!white, fill opacity =0.2]
 (0,-10) circle (7.5cm);
\filldraw[draw=black!50!white, fill=black!30!white, fill opacity=0.3]
 (-2.5,0.1) circle (1.2 cm) (2.5,0) circle (1.2 cm);
 \filldraw[draw=black!50!white, fill=black!30!white, fill opacity=0.4]
(0,-1.8) circle (1.2cm);

\path (-2,-3.3) node[] {$B(e;nR)$};
\path (-2.7,-2.4) node[above] {$B(e;nR+11M)$};

\draw[black, very thin, dashed] plot [smooth] coordinates {(-2.5,0.5) (-1.5,-0.6) (-0.4,-1.5) (0.4,-1.5) (1.2,-0.8) (2,0.2)};
\begin{scope}
 \clip (1.2,-0.8) rectangle (2,0.2);
 {
 \draw[black, semithick] plot [smooth] coordinates {(-2.5,0.5) (-1.5,-0.6) (-0.4,-1.5) (0.4,-1.5) (1.2,-0.8) (2,0.2)};
 }
\end{scope}

\draw[black, very thin, dashed] plot [smooth] coordinates {(-2.5,0.5) (-1.5,0) (-0.4,-0.5) (0.4,-0.5) (1.2,-0.1) (2,0.2)};
\begin{scope}
 \clip (0.4,-0.5) rectangle (2,0.2);
 {
 \draw[black, semithick] plot [smooth] coordinates {(-2.5,0.5) (-1.5,0) (-0.4,-0.5) (0.4,-0.5) (1.2,-0.1) (2,0.2)};
 }
\end{scope}

\draw[black, semithick] (0,-3) -- (0.05, -2.8);
\draw[black, semithick, ->] (0,-3) -- (0, -3.2);

\path (-3.2,0.1) node[right] {$X_i$};
\path (3.2,0) node[left] {$X_j$};
\path (0,-2.4) node[above] {$X_{k'}$/$X_{k}$};
\node[vertex]
(y) at (1.1,-0.8){};
\path (y) node[below right] {$P''$};
\node[vertex]
(y') at (1,-1.3){};

\node[vertex]
(z) at (0.5,-0.5){};
\path (z) node[above left] {$P'$};
\node[vertex]
(z') at (0.3,-0.9){};

\draw[black, semithick] (y) -- (y') (z) -- (z');

\node[vertex]
(e_k) at (0,-3){};
}
\path (0, -3.2) node[below] {$e$};
\path (-1.45, -0.5) node[right] {$P$};
\end{tikzpicture}

 \caption{Avoiding bottlenecks when $P\cap B\neq\emptyset$}\label{figglue4}
\end{figure}
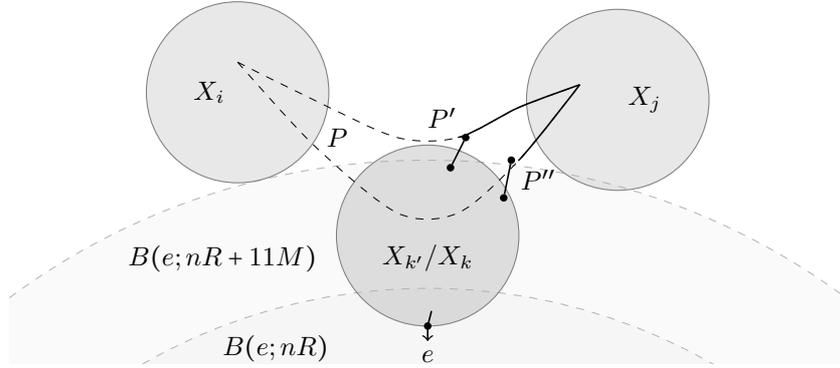
In particular, there must be some $l\in I_{j,e}$ and points $w_1,w_2\in W_{j,e}\cap X_l$ with $d(w_1,e)\geq nR+6M$, $d(w_2,e)\leq nR+M$ and such that every path from $X_j$ to $e$ meets $B_1\coloneqq B(w_1;M)$ and $B_2\coloneqq B(w_2;M)$. If $X_l=X_e$ then we set $w_2=e$. Notice that $d_X(e_l,e)\leq d_X(e_l,w_2) + d_X(w_2,e)< nR+3M$.

The following two paths from $X_k$ to $X_l$ are at distance at least $2M$, so we deduce that $k=l$.
\begin{itemize}
\item $P_1$ (avoids $B(e;nR+5M)$): take a geodesic from $w_1$ to some point in $B_1\cap\geo{g_j}$ and follow $\geo{g_j}$ to $e_j$, then join this via a path in $X_j$ to the end of $P$ contained in $X_j$, follow $P$ to $N_{4M}(X_k)$ and take any path of length at most $4M$ into $X_k$.
\item $P_2$ (contained in $B(e;nR+3M)$): take $\geo{g_l}\circ\ogeo{g_k}$.
\end{itemize}
\vspace{-3mm}
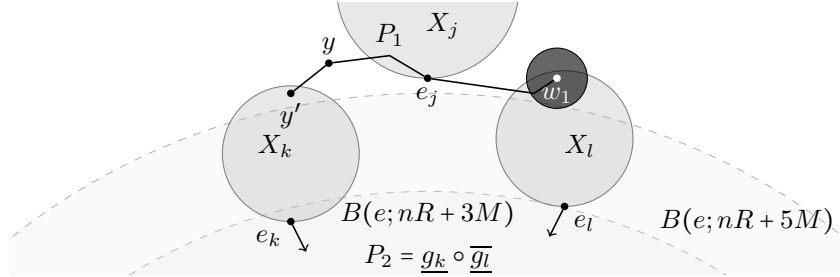
\begin{figure}[H]
 \centering
 \begin{tikzpicture}[xscale=1, yscale=1, vertex/.style={draw,fill,circle,inner sep=0.3mm}]
\clip (-5.5,-3.6) rectangle (5.5,0);
{
 \filldraw[draw=black!30!white, dashed, very thin, fill=black!10!white, fill opacity =0.2]
 (0,-10) circle (7.5cm);
  \filldraw[draw=black!30!white, dashed, very thin, fill=black!10!white, fill opacity =0.2]
 (0,-10) circle (8.8cm);
\filldraw[draw=black!50!white, fill=black!30!white, fill opacity=0.3]
 (0,0.2) circle (1.2 cm) (-1.8,-2) circle (0.9cm) (1.8,-1.8) circle (0.9cm);

\path (0,-2.8) node[] {$B(e;nR+3M)$};
\path (4.2,-3.2) node[above] {$B(e;nR+5M)$};

\node[vertex] (B) at (-1.8,-2.9) {};
\path (B) node[below left] {$e_k$};
\node[vertex] (D) at (1.8,-2.7) {};
\path (D) node[below right] {$e_l$};
\draw[black, semithick, ->] (B) -- +(0.2,-0.4);
\draw[black, semithick, ->] (D) -- +(-0.2,-0.4);
\path (0,-3.4) node[] {$P_2=\geo{g_k}\circ\ogeo{g_l}$};

\path (-2,-1.9) node[] {$X_k$};
\path (2,-1.9) node[] {$X_l$};
\path (0.2,-0.3) node[] {$X_j$};
\node[vertex]
(y) at (-1.3,-0.8){};
\path (y) node[above] {$y$};
\node[vertex]
(y') at (-1.8,-1.2){};
\path (y') node[below] {$y'$};
\node[vertex]
(e_j) at (0,-1){};
\path (e_j) node[below] {$e_j$};
\draw[black, semithick] (e_j) -- (-0.5,-0.7) -- (y) -- (y');

\filldraw[fill=black, fill opacity=0.6]
 (1.7,-1) circle (0.4 cm);
 \node[vertex, color=white]
(c) at ( 1.7, -1) {};
\path (c) node[below, color=white] {$w_1$};
\draw[black, semithick] (e_j) -- (1.4,-1.2) -- (c);
\path (-0.5,-0.7) node[above] {$P_1$};

}

\end{tikzpicture}

 \caption{Paths $P_1$ and $P_2$}\label{figsamepiece}
\end{figure}
Hence, $\geo{g_j}\cap B_1\neq\emptyset$. This completes the first part of the proof.

Now we prove that $\setcon{c_j}{j\in[i]}\subseteq N_{4M}(X_k)$. Recall that $c_j$ is the unique point on $\geo{g_j}$ at distance $nR$ from $e$.

Using the argument of the first part, we let $m_j\in \geo{g_j}\cap B_M(w_j)$, where $d(e,w_j)\geq nR-M$.
\m
If $d_X(e,w_j)\leq nR+2M$ then $d_X(m_j,e)\in(nR-2M,nR+3M)$, which implies that $d_X(w_j,c_j)\leq d_X(w_j,m_j)+d_X(m_j,c_j) < M+3M = 4M$ as required.
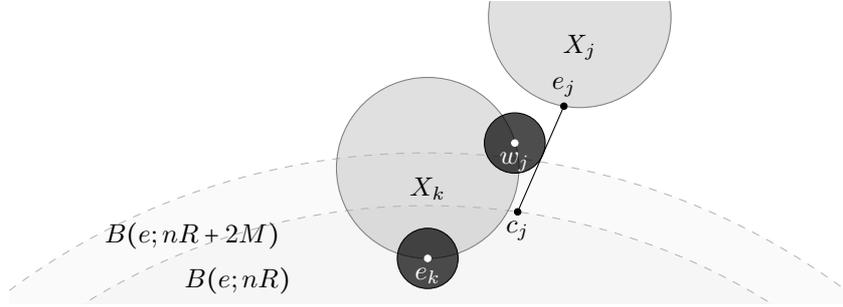
\begin{figure}[H]
 \centering
 \begin{tikzpicture}[xscale=1, yscale=1, vertex/.style={draw,fill,circle,inner sep=0.3mm}]
 \clip (-5.5,-2.8) rectangle (5.5,1.2);
{
 \filldraw[draw=black!30!white, dashed, very thin, fill=black!10!white, fill opacity =0.2]
 (0,-10) circle (9.2cm);
 \filldraw[draw=black!30!white, dashed, very thin, fill=black!10!white, fill opacity =0.2]
 (0,-10) circle (8.5cm);
\filldraw[draw=black!50!white, fill=black!30!white, fill opacity=0.4]
 (0,-1) circle (1.2 cm) (2,1) circle (1.2 cm);
 \node (A) at (2,1) {};
\node[vertex]
(B) at ([shift={(-100:1.2 cm)}]A) {};
\path (B) node[above] {$e_j$};
\node (C) at (0,-10) {};
\node[vertex]
(D) at ([shift={(82:8.5 cm)}]C) {};
\path (D) node[below] {$c_j$};
\draw[black] (B) -- (D);
\node (E) at ([shift={(83:9.4 cm)}]C) {};
\filldraw[fill=black, fill opacity=0.7]
 (E) circle (0.4 cm);
\node[vertex, color=white] at (E) {};
\path (E) node[below, color=white] {$w_j$};
\node (F) at (0,-2.2) {};
\filldraw[fill=black, fill opacity=0.7]
 (F) circle (0.4 cm);
\node[vertex, color=white] at (F) {};
\path (F) node[below, color=white] {$e_k$};
\path (-2.5,-2.5) node[] {$B(e;nR)$};
\path (-3.1,-2.2) node[above] {$B(e;nR+2M)$};
\path (2,0.6) node[] {$X_j$};
\path (0,-1) node[below] {$X_k$};
}
\end{tikzpicture}

 \caption{Using $\geo{g_j}\cap B\neq\emptyset$ when $d_X(e,w_j)\leq nR+2M$}\label{figglue1}
\end{figure}
Otherwise, $d_X(e,w_j)> nR+2M$. Every path from $X_k$ to $e$ meets $B_k\coloneqq B(e_k;M)$, and there is a path from $X_j$ to $X_k$ avoiding this ball - follow $\geo{g_j}$ from $e_j$ to $m_j$ then take a geodesic from $m_j$ to $w_j$ - so every path from $X_j$ to $e$ must meet $B_k$. Moreover, $d_X(e_k,e)\leq nR$ since $k\in I^n$. Thus, applying Lemma \ref{qconv} to the geodesic $\geo{g_j}$, we see that $c_j\in N_{4M}(X_k)$. 
\end{proof}

Given some equivalence class of pieces $[i]$, we collect all bottlenecks separating this collection of pieces from $e$ using the following definition: 
\[
W^{[i]}\coloneqq \bigcup_{j\sim i} W_{j,e} \cap \bigcup_{k'\in I^n} X_{k'}.
\]

Now, we choose a function $c: I\setminus\set{e}\to I$ which is level decreasing - $c(I_{n+1})\subseteq I^n$ - with the following properties:
\begin{itemize}
 \item if $i\sim j$ then $c(i)=c(j)$,
 \item if $i\in I_{n+1}$ and $c(i)=k$, then there exists some $i'\sim i$ and some point $w\in X_k\cap W_{i',e}$ with $d_X(w,e)\geq nR-M$ such that for all $w'\in W^{[i]}$, we have $d_X(w',e)\leq d_X(w,e)+M$.
\end{itemize}
  Note that by Lemmas \ref{highbottleneck} and \ref{glue} for each $i$ there is some $k$ satisfying the above property.
  
  We give a useful criterion for determining the value of $c$ at the end of the section but note here that in Section \ref{pfthm1} we will require this more complicated definition.
  
  In general the function $c$ is not uniquely determined by these two properties, so choices must be made. Also, $W^{[i]}$ could be infinite, so $\sup(d_X(w,e))$ is not necessarily attained.

We now give the definition of a tree-graded space $\cT{X}$ associated to $X$.

The space $\cT{X}$ is defined inductively starting with a base space $\cT{X}_0=N_{4M}(X_e)$. We construct $\cT{X}_k$ from $\cT{X}_{k-1}$ by adding a copy of $N_{4M}(X_i)$ for each $i\in I_k$ and attach $e_i\in N_{4M}(X_i)$ to $c_i\in N_{4M}(X_{c(i)})$ by a geodesic of length $d_X(e_i,c_i)$. By Lemma \ref{glue}, this construction is well-defined.

Defining $\cT{X}=\bigcup_{k\in\N}\cT{X}_k$ gives a tree-graded space whose set of pieces consists of singleton sets and $\setcon{\mcT_i\!\coloneqq N_{4M}(X_i)}{i\in I}$. We denote the natural metric on $\cT{X}$ by $d_{\cT{X}}$.

The underlying tree $\mathcal{T}$ for this construction is defined to have vertex set $I$ and $ij$ is an edge if and only if $c(i)=j$ or $c(j)=i$. The simplicial graph metric on $\mathcal{T}$ is denoted by $d_{\mcT}$.

We make one important observation at this point. If $X$ is a simplicial graph, $M\in\Z$ and we choose base points $e_i$ which are vertices, then it is easy to give $\cT{X}$ the structure of a simplicial graph by dividing the (integer length) paths $e_ic_i$ into edges of length $1$.

We finish this section with a criterion which is sufficient to determine $c(i)$.

\begin{lem}\label{c(i)} \h If $\lv{i}\coloneqq n+1>\lv{j}$ and there exists some path $P$ from some $X_i$ to $X_j$ avoiding $B(e;nR+5M)$ then $c(i)=j$.
\end{lem}
\begin{proof} Firstly, $j\in I_{i,e}$. This follows exactly from the proof of Lemma \ref{glue}; we find some suitable $X_l$ with $l\in I_{i,e}\cap I^n$ and a point $w_r\in W_{i,e}\cap X_l$ with $d(e,w_r)\geq nR+4M$, then prove $j=l$.
\m
Now suppose $c(i)=j'\neq j$, so there is some $i'\sim i$ and $j'\in I_{i',e}\cap I^n$ such that $c(i')=c(i)=j'$. By definition, $X_{j'}$ contains a point $w\in W_{i',e}$ with $d_X(w,e)\geq nR+3M$ such that all paths from $X_{i'}$ to $X_e$ meet $B(w;M)$.
\m
Let $P_0$ be some path from $X_i$ to $X_{i'}$ with $P_0\cap B(e;nR+11M)\subseteq N_{4M}(X_{k})$ for some $k\in I^n$ and consider the paths $P_1,P_2$ from $X_j$ to $X_{j'}$ given below: 
\begin{itemize}
 \item $P_1$: (contained in $B(e;nR)$) concatenate $\geo{g_j}$ with $\ogeo{g_{j'}}$,
 \item $P_2$: start at $w_{r}$ and take a path of length at most $M$ to some $m_i\in\geo{g_i}$ then follow the reverse of $\geo{q_{y_i}^i}$ to $y_i$, take $P_0$ to $y_{i'}$, $\geo{q_{y_{i'}}^{i'}}$ to some $m_{i'}\in \geo{g_{i'}}\cap B(w;M)$ then take some path of length at most $M$ to $w$.
\end{itemize}
\vspace{-3mm}
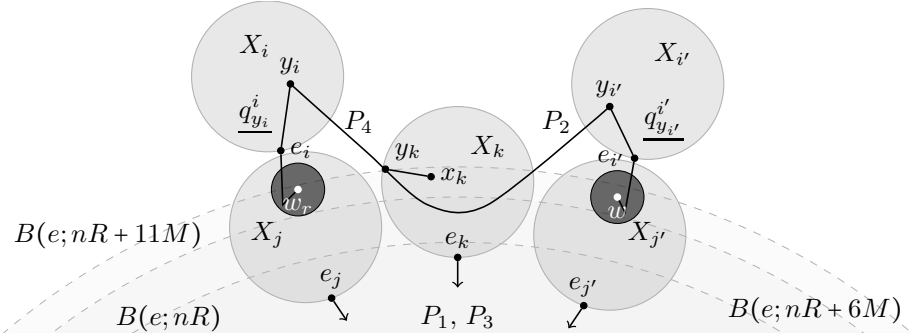
\begin{figure}[H]
 \centering
 \begin{tikzpicture}[xscale=1, yscale=1, vertex/.style={draw,fill,circle,inner sep=0.3mm}]
\clip (-6,-3) rectangle (6,1.4);
{
\filldraw[draw=black!30!white, dashed, very thin, fill=black!10!white, fill opacity =0.15]
 (0,-10) circle (9.2cm);
 \filldraw[draw=black!30!white, dashed, very thin, fill=black!10!white, fill opacity =0.15]
 (0,-10) circle (8.8cm);
 \filldraw[draw=black!30!white, dashed, very thin, fill=black!10!white, fill opacity =0.15]
 (0,-10) circle (8.2cm);
\filldraw[draw=black!30!white, fill=black!30!white, fill opacity=0.3]
 (-2.5,0.4) circle (1 cm) (2.5,0.3) circle (1 cm) (-2,-1.6) circle (1 cm) (2,-1.7) circle (1 cm) (0,-1) circle (1 cm);
 
\node (A) at (-2.5,0.4) {};
\node[vertex]
(B) at ([shift={(-80:1 cm)}]A) {};
\path (B) node[right] {$e_i$};
\path (B) +(0,0.1) node[above left] {$\geo{q_{y_i}^i}$};
\node (C) at (2.5,0.3) {};
\node[vertex]
(D) at ([shift={(-100:1 cm)}]C) {};
\path (D) node[left] {$e_{i'}$};
\path (D) +(0,0.1) node[above right] {$\geo{q_{y_{i'}}^{i'}}$};
\node (E) at (-2,-1.6) {};
\node[vertex]
(F) at ([shift={(-70:1 cm)}]E) {};
\path (F) node[above] {$e_j$};
\node[vertex]
(G) at (0,-2) {};
\path (G) node[above] {$e_{k}$};
\node (H') at (2,-1.7) {};
\node[vertex]
(H) at ([shift={(-110:1 cm)}]H') {};
\path (H) node[above] {$e_{j'}$};

\path (-3.8,-2.8) node[] {$B(e;nR)$};
\path (4.7,-2.7) node[] {$B(e;nR+6M)$};
\path (-4.6,-1.7) node[] {$B(e;nR+11M)$};

\filldraw[fill=black, fill opacity=0.55]
 ( -2.1, -1.1) circle (0.35 cm);
\filldraw[fill=black, fill opacity=0.55]
 ( 2.1, -1.2) circle (0.35 cm);

\draw[black, semithick] plot [smooth] coordinates {(-2.2,0.3) (-1.2,-0.6) (-0.4,-1.3) (0.4,-1.3) (2,0)};
\draw[black, semithick] (-2.2,0.3) -- (B) -- (-2.3, -1.3) -- (-2.1,-1.1) (2,0) -- (D) -- (2.2,-1.4) -- (2.1,-1.2);

\node[vertex, color=white]
(K) at ( -2.1, -1.1) {};
\path (K) node[below, color=white] {$w_r$};
 \node[vertex, color=white]
(L) at ( 2.1, -1.2) {};
\path (L) node[below, color=white] {$w$};

\node[vertex]
(I) at (-2.2,0.3) {};
\path (I) node[above] {$y_i$};
\node[vertex]
(J) at (2,0) {};
\path (J) node[above] {$y_{i'}$};
\node[vertex]
(K) at (-0.95,-0.83) {};
\path (K) node[above right] {$y_{k}$};
\node[vertex]
(L) at ([shift={(0.6,-0.1)}]K) {};
\path (L) node[right] {$x_{k}$};

\path (-3,0.8) node[right] {$X_i$};
\path (3.2,0.7) node[left] {$X_{i'}$};
\path (0.4,-0.8) node[above] {$X_{k}$};
\path (2.5,-1.7) node[] {$X_{j'}$};
\path (-2.5,-1.7) node[] {$X_{j}$};
\path (1.3,-0.2) node[] {$P_2$};
\path (-1.3,-0.2) node[] {$P_4$};
\path (0,-2.8) node[] {$P_1$, $P_3$};

\draw[black, semithick] (K) -- (L);
\draw[black, semithick, ->] (F) -- +(0.2,-0.3);
\draw[black, semithick, ->] (G) -- +(0,-0.4);
\draw[black, semithick, ->] (H) -- +(-0.2,-0.3);
}
\end{tikzpicture}

 \caption{Paths $P_1$ and $P_2$}\label{figtwistpath}
\end{figure}
  These paths are at distance at least $2M$ - contradicting \tu{(RBP)} - unless there is some $p\in P_0$ with $ d(p,P_1)<2M$. It follows that $p\in B(e;nR+M)$, and therefore $p\in N_{4M}(X_{k})$.
\m
In this situation we prove $j=k=j'$, we present only the first of these, the second follows using the same method. To do this we present two paths $P_3$ and $P_4$ from $X_j$ to $X_{k}$ at distance at least $2M$ (cf. Figure \ref{figtwistpath}).
\begin{itemize}
 \item $P_3$: (contained in $B(e;nR)$) concatenate $\geo{g_j}$ with $\ogeo{g_{k}}$.
 \item $P_4$: (avoids $B(e;nR+2M)$) follow $P_2$ from $w_{r}$ to a suitable point $y_{k}\in P_0\cap N_{4M}(X_{k})$ then take any path of length at most $4M$ to some point $x_{k}\in X_{k}$.
\end{itemize}
This completes the proof. \end{proof}

\section{A quasi-isometry from $\cT{X}$ to $X$}\label{pfthm1}
  Here we show that the natural collapse $\phi\!:\cT{X}\to X$ which maps each $\mcT_i$ onto $N_{4M}(X_i)$ in the obvious way defines a quasi-isometry.

From the construction it follows immediately that $\phi$ is $1$-Lipschitz and surjective.

We denote by $e'_i$ and $c'_i$ the unique points in $\cT{X}$ contained in $\phi^{-1}(e_i)\cap\mcT_i$ and  $\phi^{-1}(c_i)\cap\mathcal{T}_{c(i)}$ respectively.

To prove the other inequality we take any two points $x\in \mcT_i$ and $y\in \mathcal{T}_j$ and write the $\mathcal{T}$-geodesic between $i$ and $j$ as
\[
i=i_0,i_1,\dots,i_a=l=j_b,j_{b-1},\dots, j_0=j,
\]
where $l$ is the unique piece along this geodesic of minimal level.

Without loss of generality we may assume $d_X(e_i,e)\geq d_X(e_j,e)$.

We firstly deal with the case where at least one of $a,b$ is $0$. By our above assumption, it must be the case that $b=0$. To achieve this we present a base case (Lemma \ref{a1b0}) and then apply an inductive process on $a$ (Lemma \ref{b0}).

\begin{lem}\label{a1b0} \h Suppose in the above situation $a\leq 1$ and $b=0$, then
\[
d_{\cT{X}}(x,y) \leq d_{X}(\phi(x),\phi(y)) + 2R + 32M.
\]
\end{lem}
\begin{proof} If $a=0$ then $i=j$ and the result is obvious as $X_i$ is $4M$ quasi-convex. For $a=1$, $\lv{j}<\lv{i}$ so by Lemma \ref{I^nlv}, any path from $\phi(x)$ to $\phi(y)$ meets $B(e_i;8M)$. Hence, $d_X(\phi(x),\phi(y))\geq d_X(\phi(x),e_i) + d_X(e_i,\phi(y)) - 16M$. Moreover,
\[
\begin{array}{rcl}
d_{\cT{X}}(x,y)
&
=
&
d_{\cT{X}}(x,e'_i) + d_{\cT{X}}(e'_i,c'_i) + d_{\cT{X}}(c'_i,y)
\m
&
\leq 
&
(d_X(\phi(x),e_i) + 8M) +d_X(e_i,c_i) + (d_X(c_i,\phi(y)) + 8M)
\m
&
\leq
&
d_X(\phi(x),e_i) + d_X(e_i,\phi(y)) + 16M + 2R.
\end{array}
\]
The result follows by combining the two inequalities. \end{proof}

Our first inductive step completes the proof in the case $b=0$.

\begin{lem}\label{b0} \h Suppose $a\geq 2$ and $b=0$. Then
\[
d_{\cT{X}}(x,y) \leq d_{X}(\phi(x),\phi(y)) + 2R + 72Ma + 16M.
\]
\end{lem}
\begin{proof} Note that by construction there is some $i'\sim i$ such that $c(i)\in I_{i',e}$, and the $w\in X_{c(i)}\cap W_{i',e}$ with $d(w,e)$ maximal satisfies $d(w,e)\geq nR-M$. Set $B=B(w;5M)$.
 
Firstly, we prove that every geodesic $\geo{g}\in\setg{\phi(x)}{\phi(y)}$ meets $N_{5M}(X_{c(i)})$.

Suppose that some geodesic $\geo{g}\in\setg{\phi(x)}{\phi(y)}$ avoids $B$, then we consider all ways of using the relation $i\sim i'$ and the geodesic $\geo{g_j}$ to extend $\geo{g}$ to a path $P$ from $X_{i'}$ to $e$ and deduce that $d(w,e)\geq nR+6M$.
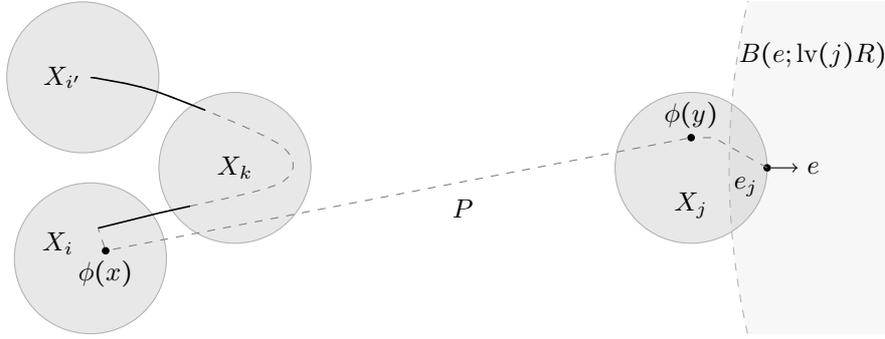
\begin{figure}[H]
 \centering
 \begin{tikzpicture}[xscale=1, yscale=1, vertex/.style={draw,fill,circle,inner sep=0.3mm}]
\clip (-5.4,-2.2) rectangle (6.6,2.2);
{

\filldraw[draw=black!30!white, dashed, fill=black!10!white, fill opacity=0.3]
 (14.5,0) circle (10 cm);

\filldraw[draw=black!30!white, fill=black!30!white, fill opacity=0.3]
 (-4,1.2) circle (1 cm) (-3.9,-1.2) circle (1 cm) (-2,0) circle (1 cm) (4,0) circle (1 cm);

\node[vertex]
(F) at (5,0) {};
\path (F) node[below left] {$e_j$};
\path (5.4,0) node[right] {$e$};

\draw[black!50!white, thin, dashed] plot [smooth] coordinates {(-3.9,1.2) (-3,1) (-1.4,0.3) (-1.4,-0.2) (-3,-0.6) (-3.8,-0.8)};

\node[vertex]
(I) at (4,0.4) {};
\path (I) node[above] {$\phi(y)$};
\node[vertex]
(J) at (-3.7,-1.1) {};
\path (J) node[below] {$\phi(x)$};

\draw[black!50!white, dashed] (5,0) -- (4.3,0.4) -- (I) -- (J) -- (-3.8,-0.8);

\path (-4,-1) node[left] {$X_i$};
\path (-3.9,1.2) node[left] {$X_{i'}$};
\path (4,-0.5) node[] {$X_{j}$};
\path (-2,0) node[] {$X_k$};
\path (1,-0.3) node[below] {$P$};

\begin{scope}
 \clip (-3.8,-0.8) rectangle (-2.6, -0.5);
 {
 \draw[black, semithick] plot [smooth] coordinates {(-3.9,1.2) (-3,1) (-1.4,0.3) (-1.4,-0.2) (-3,-0.6) (-3.8,-0.8)};
 }
\end{scope}

\begin{scope}
 \clip (-4,1.2) rectangle (-2.4, 0.5);
 {
 \draw[black, semithick] plot [smooth] coordinates {(-3.9,1.2) (-3,1) (-1.4,0.3) (-1.4,-0.2) (-3,-0.6) (-3.8,-0.8)};
 }
\end{scope}

\draw[black, thin, ->] (F) -- +(0.4,0);

\path (5.6,1.5) node[] {$B(e;\lv{j}R)$};
}
\end{tikzpicture}
 \caption{A path $P$ and possible intersection with $B(w;M)$}\label{figavoidws}
\end{figure}
  Therefore, $c(i)\in I_{i,e}$ using the proof of Lemma \ref{glue}. In particular, this means that every path from $X_i$ to $e$ meets $B(w';M)$ for some $w'\in X_{c(i)}$ with $d(e,w')\geq nR+5M$. Suppose $\geo{g}\cap B(w';5M)=\emptyset$, then we may extend $\geo{g}$ to a path from $X_i$ to $e$ avoiding $B(w';M)$, which is a contradiction. This completes the proof of the claim. 

Every geodesic from $\phi(x)$ to $\phi(y)$ meets $N_{5M}(X_{c(i)})$, so they also meet $B'=B(e_{c(i)};9M)$ by Lemma \ref{I^nlv}. We now build a $62M$-slack geodesic $\geo{q}$ from $\phi(x)$ to $\phi(y)$ which contains $c_i$.

Follow $q^{\phi(x)}_i$ (cf. Lemma \ref{10Mslack}) from $\phi(x)$ to some point $z$ after $c_i$ which is contained in $B(e_{c(i)};8M)$, (by Lemmas \ref{glue} and \ref{I^nlv}, any path from $c_i$ to $e$ meets $B(e_{c(i)};8M)$); take a geodesic from $z$ to some point $z'\in\geo{g}\cap B'$ and then follow $\geo{g}$ to $\phi(y)$.

Now $d(\phi(x),\phi(y)) \geq d(\phi(x),z) + d(z,z') + d(z',\phi(y)) - 52M$, so using the fact that $q^{\phi(x)}_i$ is $10M$-slack, we deduce that $\geo{q}$ is a $62M$-slack geodesic.
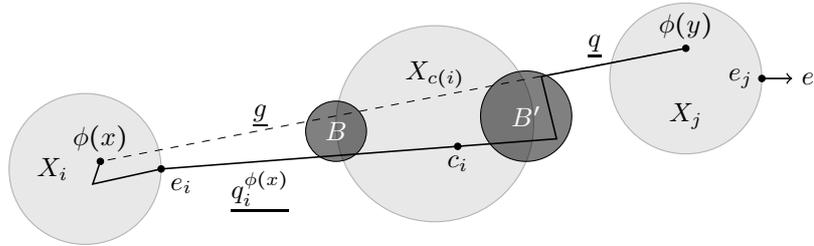
\begin{figure}[H]
 \centering
 \begin{tikzpicture}[xscale=1, yscale=1, vertex/.style={draw,fill,circle,inner sep=0.3mm}]
\clip (-6,-2.2) rectangle (6,1);
{
\filldraw[draw=black!30!white, fill=black!30!white, fill opacity=0.3]
(-3.9,-1.2) circle (1 cm) (0.7,-0.6) circle (1.3 cm) (4,0) circle (1 cm);

\node[vertex]
(F) at (5,0) {};
\path (F) node[left] {$e_j$};
\path (5.4,0) node[right] {$e$};

\filldraw[fill=black, fill opacity=0.5]
 ( -0.6, -0.7) circle (0.4 cm);
 \filldraw[fill=black, fill opacity=0.5]
 ( 1.9, -0.5) circle (0.6 cm);

\path (1.9,-0.5) node[color=white] {$B'$};

\path (-0.6,-0.7) node[color=white] {$B$};

\node[vertex]
(I) at (4,0.4) {};
\path (I) node[above] {$\phi(y)$};
\node[vertex]
(J) at (-3.7,-1.1) {};
\path (J) node[above] {$\phi(x)$};
\node[vertex]
(K) at (1,-0.9) {};
\path (K) node[below] {$c_i$};

\draw[black, very thin, dashed] (I) -- (J);
\draw[black, semithick] (J) -- +(-0.1,-0.3) -- (-2.9,-1.2) -- (1,-0.9) -- (2.3, -0.8) -- (2.1, 0.025) -- (I);

\path (-4,-1.2) node[left] {$X_i$};
\path (0.7,0.05) node[] {$X_{c(i)}$};
\path (4,-0.5) node[] {$X_{j}$};
\path (-1.6,-0.5) node[] {$\geo{g}$};
\path (-1.6,-1.1) node[below] {$\geo{q_i^{\phi(x)}}$};
\path (2.8,0.15) node[above] {$\geo{q}$};
\node[vertex]
(L) at (-2.9,-1.2) {};
\path (L) node[below right] {$e_i$};

\draw[black, semithick, ->] (F) -- +(0.4,0);
}
\end{tikzpicture}

 \caption{The $62M$-slack geodesic $\geo{q}$}\label{fig62Mslack}
\end{figure}
  Since $\geo{q}$ meets $c_i$,
\[
d_{X}(\phi(x),\phi(y)) \geq d_X(\phi(x),c_i)+ d_X(c_i,\phi(y))-62M.
\]
We recall that by the inductive hypothesis,
\[
d_{\cT{X}}(c'_i,y) \leq d_{X}(c_i,\phi(y)) + 2R + 72M(a-1) + 16M.
\]
Finally, by Lemma \ref{10Mslack}, $d_{\cT{X}}(x,c'_i) = d_X(\phi(x),e_i)+d_X(e_i,c_i)\leq d_X(\phi(x),c_i)+10M$, so combining these we see that
\[
\begin{array}{rcl}
d_{\cT{X}}(x,y)
&
=
&
d_{\cT{X}}(x,c'_i) + d_{\cT{X}}(c'_i,y)
\m
& 
\leq 
&
d_X(\phi(x),\phi(y)) + 2R + 72M(a-1) + 62M + 10M + 16M 
\m
&
= 
&
d_X(\phi(x),\phi(y)) + 2R + 72Ma + 16M.
\end{array}
\]
\end{proof}

Now we come to the case $b\geq 1$. Again we start with a base case before progressing to the general result.

\begin{lem}\label{ab1} \h Suppose $a=b=1$. Then
\[
d_{\cT{X}}(x,y)\leq d_{X}(\phi(x),\phi(y))+7R+70M.
\]
\end{lem}
\begin{proof} Recall that $l=c(i)=c(j)$. Without loss of generality we assume $d(e_i,e)\geq d(e_j,e)$, so in particular, $n\coloneqq\lv{i}\geq m\coloneqq \lv{j}$. By Lemma \ref{I^nlv}, every path from $\phi(x)$ to $N_{4M}(X_j)$ passes through $B(e_i;8M)$. If some geodesic in $\setg{\phi(x)}{\phi(y)}$ meets $B(e_j;14M)$, then 
\[
\begin{array}{rcl}
d_X(\phi(x),\phi(y))
&
\geq 
&
d_X(\phi(x),e_j)+d_X(e_j,\phi(y))-28M 
\m
&
\geq 
&
d_X(\phi(x),e_i)+ d_X(e_i,e_j) - 14M + d_X(e_j,\phi(y)) - 28M.
\end{array}
\]
Combining these bounds we see that
\[
\begin{array}{rcl}
d_{\cT{X}}(x,y) 
&
=
&
d_{\cT{X}}(x,e'_i) + d_{\cT{X}}(e'_i,c'_i) + d_{\cT{X}}(c'_i,c'_j) +
\m
&
&
d_{\cT{X}}(c'_j,e'_j) + d_{\cT{X}}(e'_j,y)
\m
&
\leq
&
d_X(\phi(x),e_i) + d_X(c_i,c_j) + d_X(e_j,\phi(y)) + 2R + 28M
\m
&
\leq
&
d_X(\phi(x),e_i) + d_X(e_i,e_j) + d_X(e_j,\phi(y)) + 4R + 28M
\m
&
\leq
&
d_X(\phi(x),\phi(y)) + 4R + 70M.
\end{array}
\]
Now suppose all geodesics avoid $B(e_j;14M)$. By Lemma \ref{I^nlv} we know that geodesics must also avoid $\bigcup_{k\in I^m} N_{10M}(X_k)$, so, in particular they avoid the set $N_{6M}(\geo{g_i^c})$ where we define $\geo{g_i^c}$ to be the restriction of $\geo{g_i}$ to a geodesic in $\setg{c_i}{e}$. Moreover, all geodesics must also avoid $N_{6M}(\geo{g_j})$ otherwise one can find a path from $X_j$ to $e$ avoiding $B(e_j;M)$.

Hence, the bottleneck $w=w_0\in W_{j,i}$ lying in $X_j$ must be within distance $M$ of some point of $\geo{g_i}\setminus B(e;nR+6M)$. In particular there is a path from $X_i$ to $X_j$ avoiding $B(e;nR+5M)$.
\begin{figure}[H]
 \centering
 \begin{tikzpicture}[xscale=1, yscale=1, vertex/.style={draw,fill,circle,inner sep=0.3mm}]
\clip (-6,-2.2) rectangle (6,1);
{
\node[vertex]
(I) at (5,0) {};
\node[vertex]
(N) at (-1.5,-0.99) {};

\filldraw[draw=black!30!white, dashed, very thin, fill=black!20!white, fill opacity =0.2]
 (I) circle (6.4cm);

\fill[black!30!white]
 (N) circle (0.4cm) (I) circle (0.4cm);
\fill[black!30!white] 
 (5,-0.4) -- (-1.5,-1.39) -- (-1.5,-0.59) -- (5,0.4) -- (5,-0.4);
 \fill[black!30!white]
 (-1.1,0) circle (0.4cm);
\fill[black!30!white] 
 (-1.1,-0.4) -- (-1.1,0.4) -- (5,0.4) -- (5,-0.4) -- (-1.1,-0.4);

\filldraw[draw=black!30!white, fill=black!30!white, fill opacity=0.3]
 (-2.1,0) circle (1 cm) (-3.9,-1.2) circle (1 cm);

\filldraw[fill=black, fill opacity=0.7]
 (-2.1,-1) circle (0.4 cm);
\node[vertex, color=white]
(c) at ( -2.1, -1) {};
\path (c) node[above, color=white] {$w$};

\node[vertex]
(I) at (5,0) {};
\path (I) node[above] {$e$};
\node[vertex]
(J) at (-3.7,-1.1) {};
\path (J) node[left] {$\phi(x)$};
\node[vertex]
(J') at (-3.9,1.2) {};
\path (J') node[above] {$\phi(y)$};
\node[vertex]
(K) at (-1.1,0) {};
\path (K) node[below] {$e_j$};
\node[vertex]
(L) at (-2.9,-1.2) {};
\path (L) node[below right] {$e_i$};
\node[vertex]
(M) at (-0.65,0) {};
\path (M) node[below] {$c_j$};
\node[vertex]
(N) at (-1.5,-0.99) {};
\path (N) node[below right] {$c_i$};

\path (0.3,-2) node[] {$B(e;nR+7M)$};

\path (-3.9,-0.6) node[] {$X_i$};
\path (-2.5,0) node[] {$X_{j}$};

\draw[black, very thin, dashed] (K) -- (I) -- (L);
\draw[black, semithick] (J) -- (-3.9,-1.4) -- (L) -- (-1.25, -0.95) -- (-1.35,-0.64);
}
\end{tikzpicture}
 \caption{A path from $X_i$ to $X_j$}\label{figshortpath}
\end{figure}
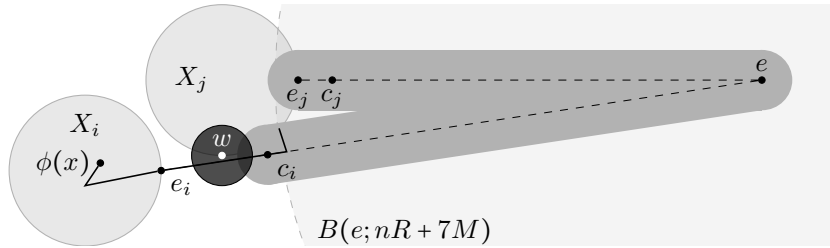
  If $n>m$ then $c(i)=j$, by Lemma \ref{c(i)}, which contradicts the assumption that $a=b=1$. However, if $n=m$ then $d_X(w_0,e_j)\leq R+3M$, because every geodesic from $w_0$ to $e$ meets $B(e_j;M)$ and $d_X(w_0,e)\leq d(e_j,e)+R+M$. Hence,
\[
 d_X(c_i,c_j)\leq d_X(c_i,w)+d_X(w,e_j)+d_X(e_j,c_j) \leq (R+M) + (R+3M) + R = 3R + 4M,
\]
while $d_X(e_i,e_j)\leq d_X(e_i,w)+d_X(w,e_j) \leq (R-7M) + (R+3M) = 2R-4M$.

Therefore,
\[
\begin{array}{rcl}
d_{\cT{X}}(x,y)
 &
\leq
 &
d_{X}(\phi(x),e_i) + R + (3R+4M) + R + d_{X}(e_j,\phi(y)) + 16M
\m
 &
\leq
 &
d_X(\phi(x),e_i) + d_{X}(e_i,\phi(y)) + (2R - 4M) + 5R + 20M
\m
 &
\leq
 &
d_X(\phi(x),\phi(y))+ 7R + 16M + 16M.
\end{array}
\]
The final step uses Lemma \ref{I^nlv}. \end{proof}

This leads to the final lemma required for the proof.

\begin{lem}\label{a+bgeq3} Suppose $a,b\geq 1$ and $d_{\mcT}(i,j)=a+b\geq 3$, then
\[
d_{\cT{X}}(x,y)\leq d_{X}(\phi(x),\phi(y))+ 9R + 80M(a+b).
\]
\end{lem}
\begin{proof} We proceed by induction on $a+b$ using the the situation $b\leq 1$ as the base case, we do not include the extra $+16M$ as we will not require the situation $a=b=0$ in our inductive step. To ease notation we set $\lv{i}\coloneqq n+1$ and $\lv{j}\coloneqq m+1$, by assumption $\lv{i},\lv{j}\geq 1$.

If some $45M$-slack geodesic from $\phi(x)$ to $\phi(y)$ meets $\set{c_i,c_j}$, (we deal with the case of $c_i$, the other case is very similar) then $d_X(\phi(x),\phi(y))\geq d_X(\phi(x),c_i)+d_X(c_i,\phi(y))-45M$.

Lemma \ref{10Mslack} gives $d_X(\phi(x),c_i)\geq d_X(\phi(x),e_i) + d_X(e_i,c_i)-10M$, while by the inductive hypothesis
\[
d_{\cT{X}}(y,c'_i)\leq d_X(\phi(y),c_i) + 9R + 80M(a+b-1).
\]
Combining these we see that
\[
\begin{array}{rcl}
d_{\cT{X}}(x,y)
&
=
&
d_{\cT{X}}(x,c'_i) + d_{\cT{X}}(c'_i,y)
\m
&
\leq
&
d_X(\phi(x),e_i)+d_X(e_i,c_i) + 8M + d_X(c_i,\phi(y))
\m
&
&
+ 9R + 80M(a+b-1)
\m
&
\leq
&
d_X(\phi(x),\phi(y))9R + 80M(a+b).
\end{array}
\]

Now suppose every $45M$-slack geodesic from $\phi(x)$ to $\phi(y)$ avoids $\set{c_i,c_j}$, then every geodesic in $\setg{\phi(x)}{\phi(y)}$ misses $N_{15M}(\geo{g_i^c}\cup\geo{g_j^c})$, where $\geo{g_k^c}$ is the restriction of $\geo{g_k}$ to a geodesic in $\setg{c_k}{e}$. If this is not the case then it is easy to find a suitable slack geodesic $\geo{q}$ which hits either $c_i$ or $c_j$.
\begin{figure}[H]
 \centering
 \begin{tikzpicture}[xscale=1, yscale=1, vertex/.style={draw,fill,circle,inner sep=0.3mm}]
\clip (-6,-2.2) rectangle (6,2.2);
{
\node[vertex]
(I) at (5,0) {};
\node[vertex]
(N) at (-1.5,-0.99) {};

\fill[black!30!white]
 (N) circle (0.4cm) (I) circle (0.4cm);
\fill[black!30!white] 
 (5,-0.4) -- (-1.5,-1.39) -- (-1.5,-0.59) -- (5,0.4) -- (5,-0.4);

\filldraw[draw=black!30!white, fill=black!30!white, fill opacity=0.3]
 (-4,1.2) circle (1 cm) (-3.9,-1.2) circle (1 cm);

\draw[black, very thin, dashed] plot [smooth] coordinates {(-3.9,1.2) (-3,1) (-1.4,0.3) (-1.4,-0.7) (-3,-0.7) (-3.7,-1.1)};
\begin{scope}
\clip (-3.9,-0.65) rectangle (-1.2, 1.2);
{
\draw[black, semithick] plot [smooth] coordinates {(-3.9,1.2) (-3,1) (-1.4,0.3) (-1.4,-0.7) (-3,-0.7) (-3.7,-1.1)};
}
\end{scope}

\node[vertex]
(I) at (5,0) {};
\path (I) node[above] {$e$};
\node[vertex]
(J) at (-3.7,-1.1) {};
\path (J) node[left] {$\phi(x)$};
\node[vertex]
(J') at (-3.9,1.2) {};
\path (J') node[above] {$\phi(y)$};
\node[vertex]
(K) at (-3,1.2) {};
\path (K) node[above right] {$e_j$};
\node[vertex]
(L) at (-2.9,-1.2) {};
\path (L) node[above right] {$e_i$};
\node[vertex]
(M) at (-1.5,0.98) {};
\path (M) node[above right] {$c_j$};
\node[vertex]
(N) at (-1.5,-0.99) {};
\path (N) node[below right] {$c_i$};

\path (-3.9,-0.6) node[] {$X_i$};
\path (-3.9,0.6) node[] {$X_{j}$};
\path (-3,-0.7) node[above right] {$\geo{g}$};
\path (-1, -0.2) node[] {$\geo{q}$};

\draw[black, very thin, dashed] (K) -- (I) -- (L);
\draw[black, semithick] (J) -- (-3.9,-1.4) -- (L) -- (-1.25, -0.95) -- (-1.35,-0.64);
}
\end{tikzpicture}
 \caption{Finding slack geodesics meeting $c_i$}\label{figfindslack}
\end{figure}
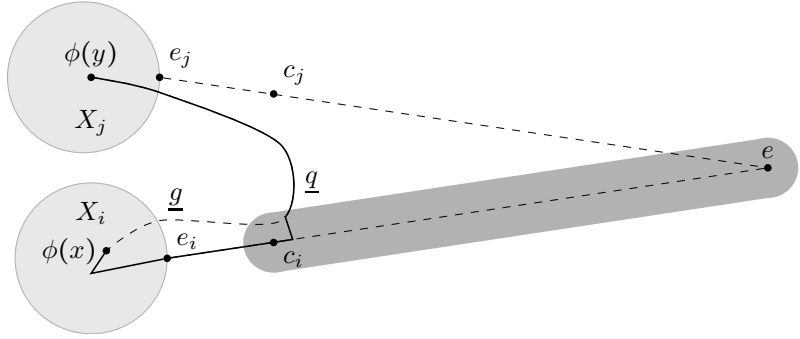
  We now have two paths from $N_{4M}(X_i)$ to $N_{4M}(X_j)$ given by $\geo{g_i}\circ\ogeo{g_j}$ and some $\geo{g}\in\setg{\phi(x)}{\phi(y)}$.

As $\geo{g}\cap N_{15M}(\geo{g_i^c}\cup\geo{g_j^c})=\emptyset$, we deduce that the collection of bottlenecks $W_{i,j}$ given by (RBP) is contained in
\[
\left(N_{M}(\geo{g_i})\setminus B(e;nR+13M)\right) \cup \left(N_{M}\left(X_j\cup \geo{g_j}\right)\setminus B(e;mR+13M)\right).
\]
We label the first of these two sets $A$ and the second one $B$. Here we are using Lemma \ref{tech} to ensure that $X_j$ is (path-)connected.

If the Hausdorff distance between $A$ and $B$ is less than $M$ then it is clear that $i\sim j$ if $\lv{i}=\lv{j}$ or $c(i)=j$, if $\lv{i}>\lv{j}$, by Lemma \ref{c(i)}. Both of these contradict the assumption that $d_{\mcT}(i,j)\geq 3$. 

Hence, there is some piece $X_k$, with $k\in I_{i,j}$ containing two bottlenecks, one in each of $A$ and $B$. We label the bottleneck point in $A$ by $w_1$ and the one in $B$ by $w_2$.

From here on we split into a number of cases depending on the relationship between $\lv{i}$, $\lv{j}$ and $\lv{k}$.
\m
{\bf Case 1:} $\lv{i}=\lv{j}$ \h It follows immediately from the above that there is a path $P$ from $X_i$ to $X_j$ with $P\cap B(e;nR+11M)\subseteq N_{4M}(X_k)$. If $k\not\in I^n$ then we obtain a path which is disjoint from $B(e;nR+11M)$, so $i\sim j$. This contradicts the assumption that $d_{\mcT}(i,j)\geq 3$.

From now on we assume $\lv{i}>\lv{j}$.
\m
{\bf Case 2:} $\lv{k}>\lv{i}$ \h As $w_1,w_2\in X_k$, we see that $d_X(w_1,e),\ d_X(w_2,e)\geq (n+1)R-M$. Hence there is a path from $X_i$ to $X_j$ avoiding $B(e;(n+1)R-2M)$. Thus, $c(i)=j$ by Lemma \ref{c(i)}.
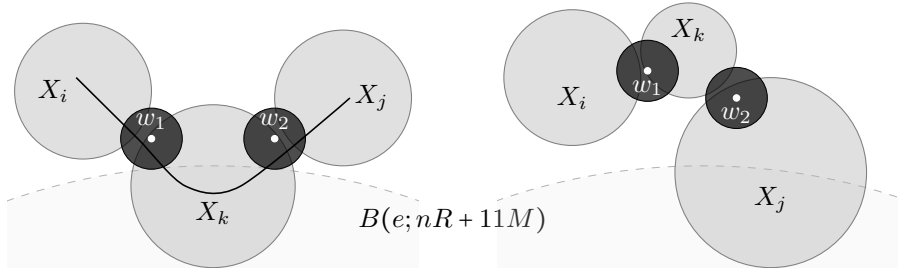
\begin{figure}[H]
 \centering
 \mbox{
\begin{tikzpicture}[xscale=0.9, yscale=0.9, vertex/.style={draw,fill,circle,inner sep=0.3mm}]
\begin{scope}
\clip (-3,-2.5) rectangle (3,1.4);
{
 \filldraw[draw=black!30!white, dashed, very thin, fill=black!10!white, fill opacity =0.2]
 (0,-10) circle (9cm);
\filldraw[draw=black!50!white, fill=black!30!white, fill opacity=0.4]
 (-1.9,0.1) circle (1 cm) (1.9,0) circle (1 cm) (0,-1.3) circle (1.2cm);

\filldraw[fill=black, fill opacity=0.7]
 (-0.9,-0.6) circle (0.45 cm);
 \node[vertex, color=white]
(c) at (-0.9,-0.6) {};
\path (c) node[above, color=white] {$w_1$};

\filldraw[fill=black, fill opacity=0.7]
 (0.9,-0.6) circle (0.45 cm);
 \node[vertex, color=white]
(d) at (0.9,-0.6) {};
\path (d) node[above, color=white] {$w_2$};

\draw[black, semithick] plot [smooth] coordinates {(-2,0.3) (-1.1,-0.6) (-0.4,-1.3) (0.4,-1.3) (2,0)};

\path (-2.7,0.1) node[right] {$X_i$};
\path (2.7,0) node[left] {$X_j$};
\path (0,-2) node[above] {$X_{k}$};
}
\end{scope}
\path (3.5,-1.8) node[] {$B(e;nR+11M)$};
\end{tikzpicture}

\hspace{-10mm}

\begin{tikzpicture}[xscale=0.9, yscale=0.9, vertex/.style={draw,fill,circle,inner sep=0.3mm}]
 \clip (-3,-2.5) rectangle (3,1.4);
{
 \filldraw[draw=black!30!white, dashed, very thin, fill=black!10!white, fill opacity =0.2]
 (0,-10) circle (9cm);
\filldraw[draw=black!50!white, fill=black!30!white, fill opacity=0.4]
 (-1.9,0.3) circle (1 cm) (-0.2,0.7) circle (0.7 cm) (1,-1.1) circle (1.4cm);

\filldraw[fill=black, fill opacity=0.7]
 (-0.8,0.4) circle (0.45 cm);
 \node[vertex, color=white]
(c) at (-0.8,0.4) {};
\path (c) node[below, color=white] {$w_1$};

\filldraw[fill=black, fill opacity=0.7]
 (0.5,0) circle (0.45 cm);
 \node[vertex, color=white]
(d) at (0.5,0) {};
\path (d) node[below, color=white] {$w_2$};

\path (-1.9,0.3) node[below] {$X_i$};
\path (1,-1.5) node[] {$X_j$};
\path (-0.2,0.7) node[above] {$X_{k}$};
}
\end{tikzpicture}
}
 \caption{Cases $1$ (left) and $2$ (right)}\label{figcase12}
\end{figure}

\noindent  {\bf Case 3:} $\lv{k}=\lv{i}$ \h Here we prove that either $c(i)=j$ or contradict the assumption that no $45M$-slack geodesic from $\phi(x)$ to $\phi(y)$ meets $c_i$.

The fact that $i\sim k$ is immediate from the location of bottleneck $w_1$.
 
If $d_X(e_k,c_k)\geq 9M$ then there is a path from $X_k$ to $X_j$ (via $w_2$) avoiding $B(e;nR+6M)$, so $c(k)=j$ by Lemma \ref{c(i)}. Hence, $c(i)=j$.

Now suppose $d_X(e_k,c_k)< 9M$, then $\geo{g_i}\cap B(e_k;M)\neq\emptyset$ as otherwise we would obtain (via $w_1$ and $\geo{g_i}$) a path from $X_k$ to $e$ avoiding $B(e_k;M)$, which contradicts (RBP). Notice that here we have used the fact that $d_X(w_1,e_k)\geq d_X(w_1,e)-d_X(e_k,e) \geq 11M-9M \geq 2M$.
 
Let $m_i\in \geo{g_i}\cap B(e_k;M)$. Then,
\[
d_X(e_k,c_i)\leq d_X(e_k,m_i)+d_X(m_i,c_i) < M + 9M + M = 11M.
\]
As every path from $\phi(x)$ to $\phi(y)$ meets $B(w_1;5M)$ it also meets $B(e_k;9M)$ by Lemma \ref{10Mslack}. Thus every such path meets $B(c_i;20M)$. In particular, there is some $40M$-slack geodesic from $\phi(x)$ to $\phi(y)$ which meets $c_i$, contradicting the initial assumption.
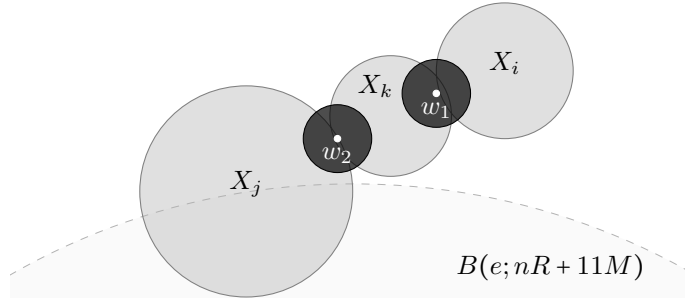
\begin{figure}[H]
 \centering
 \begin{tikzpicture}[xscale=1, yscale=1, vertex/.style={draw,fill,circle,inner sep=0.3mm}]

 \clip (-4.5,-2.5) rectangle (4.5,1.4);
{
 \filldraw[draw=black!30!white, dashed, very thin, fill=black!10!white, fill opacity =0.2]
 (0,-10) circle (9cm);
\filldraw[draw=black!50!white, fill=black!30!white, fill opacity=0.4]
 (2,0.5) circle (0.9 cm) (0.5,-0.1) circle (0.8 cm) (-1.4,-1.1) circle (1.4cm);

\filldraw[fill=black, fill opacity=0.7]
 (1.1,0.2) circle (0.45 cm);
 \node[vertex, color=white]
(c) at (1.1,0.2) {};
\path (c) node[below, color=white] {$w_1$};

\filldraw[fill=black, fill opacity=0.7]
 (-0.2,-0.4) circle (0.45 cm);
 \node[vertex, color=white]
(d) at (-0.2,-0.4) {};
\path (d) node[below, color=white] {$w_2$};

\path (2,0.6) node[] {$X_i$};
\path (-1.4,-1) node[] {$X_j$};
\path (0.3,0.3) node[] {$X_{k}$};
}
\path (2.6,-2.1) node[] {$B(e;nR+11M)$};
\end{tikzpicture}

 \caption{Case $3$: $c(i)=j$}\label{figcase3}
\end{figure}
 
From here on we assume $\lv{i}>\lv{k}$, from this and the location of the bottleneck $w_1$ we know that $c(i)=k$.
\m
{\bf Case 4:} $\lv{k}>\lv{j}$ \h As in case $3$ we find a $45M$-slack geodesic meeting $c_i$.

Immediately we see that $d_X(w_2,e_k)\leq 2M$ as the bottleneck must cut the path $\geo{g_k}\circ\ogeo{g_j}$. But as every path from $\phi(x)$ to $\phi(y)$ meets $N_{5M}(X_k)$ we see that such paths meet $B(e_k;9M)$ by Lemma \ref{10Mslack}.

Fix some $\geo{g}\in\setg{\phi(x)}{\phi(y)}$. We obtain a $45M$-slack geodesic $\geo{q}$ from $\phi(x)$ to $\phi(y)$ passing through $c_i$ by following $\geo{q^{\phi(x)}_i}$ to a point $m_i\in B(e_k;9M)\cap\geo{g}$ - if this restriction of $\geo{q_i^{\phi(x)}}$ does not include $c_i$ we include a diversion of length at most $18M$ along $\geo{q_i^{\phi(x)}}$ to $c_i$ and then back again - then follow $\geo{g}$ to $\phi(y)$.

As every path meets $B(e_k;9M)$,
\[
\begin{array}{rcl}
d_X(\phi(x),\phi(y))
&
\geq 
&
d_X(\phi(x),e_k) + d_X(e_k,\phi(y)) - 18M
\m
&
\geq
&
(d_X(\phi(x),m_i) - d_X(m_i,e_k)) + d_X(e_k,\phi(y)) - 18M
\m
&
\geq
&
l(\geo{q}) - 18M - 9M - 18M.
\end{array}
\]
(The first $-18M$ comes from the possible detour to $c_i$.) This contradicts the assumption that there is no $45M$ slack geodesic from $\phi(x)$ to $\phi(y)$ passing through $c_i$. (cf. Figure \ref{figfindslack}.)
\m
{\bf Case 5:} $\lv{j}>\lv{k}$ \h In this situation we prove that $c(i)=c(j)=k$ contradicting the assumption that $d_{\mcT}(i,j)\geq 3$.

We already know that $c(i)=k$. It is immediate from the location of $w_2$ that $d_X(e,w_2)\geq mR+10M$, so there is a path from $X_j$ to $X_k$ avoiding $B(e;mR+7M)$ and we apply Lemma \ref{c(i)} to deduce that $c(j)=k$.
\m
{\bf Case 6:} $\lv{j}=\lv{k}$ \h Here $c(i)=k\sim j$, so $d_{\mcT}(i,j)=m+n=3$. We deal with this case directly.

Firstly, $j\sim k$, as the bottleneck between $X_j\cup\geo{g_j}$ and $X_k$ yields a path from $X_j$ to $X_k$ avoiding $B(e;nR+11M)$. To avoid contradicting (RBP) for paths between $X_j$ and $X_k$ it follows that $w_2\in B(e_j;3M)\cup B(e_k;3M)$. If this is not the case then the path of length $M$ from $w_2$ to $X_j$ and $\geo{g_k}\circ\ogeo{g_j}$ are at  distance at least $2M$.

If $w_2\in B(e_j;3M)$, then 
\[
d_X(e_j,e_k)\leq d_X(e_j,w_2) + d_X(w_2,e_k) \leq 3M + (R + 5M)
\]
and if $w_2\in B(e_k;3M)$, then 
\[
d_X(e_j,e_k)\leq d_X(e_j,w_2) + d_X(w_2,e_k) \leq (R + 7M) + 3M.
\]
Here we are using the fact that any geodesic from $w_2$ to $e$ meets $B(e_k;M)$ or $B(e_j;2M)$. In either situation, $d_X(c_j,c_k)\leq 3R+10M$. 

Hence, as any path from $\phi(x)$ to $\phi(y)$ meets $B(e_k;5M)$ or $B(e_j;5M)$, by Lemma \ref{I^nlv}, 
\[
d_X(\phi(x),\phi(y))\geq d_X(\phi(x),e_k)+d_X(e_j,\phi(y)) - 2(R + 10M) - 10M.
\]
 Using Lemma \ref{a1b0} we see that
\[
 d_{\cT{X}}(x,e'_k) \leq d_{X}(\phi(x),e_k) + 2R + 40M
\]
Then as $d_{\cT{X}}(x,y) \leq d_{\cT{X}}(x,e'_k)+d_{\cT{X}}(c'_k,c'_j)+d_{\cT{X}}(e'_j,y) + 2R$,
\[
\begin{array}{rcl}
d_{\cT{X}}(x,y) 
 &
\leq 
 &
(d_{X}(\phi(x),e_k) + 2R + 40M) + (3R + 10M + 16M) + 
\m
 &
 &
(d_{X}(e_j,\phi(y)) + 8M) + 2R
\m
 &
\leq
 &
d_X(\phi(x),\phi(y)) + 9R + 104M.
\end{array}
\]
\end{proof}

We are now ready to prove the main theorem.

\subsection*{Proof of Theorem \ref{qTG}} 
The easier implication follows from Lemma \ref{RBPqi} and Proposition \ref{exTG}. From Lemmas \ref{a1b0}, \ref{b0}, \ref{ab1} and \ref{a+bgeq3} we know that for all $x\in\mcT_i$, $y\in\mcT_j$
\[
d_{\cT{X}}(x,y)\leq d_{X}(\phi(x),\phi(y)) + 9R + 80Md_{\mcT}(i,j) + 16M.
\]
Now, $d_{\cT{X}}(x,y)\geq R(\max\set{d_\mathcal{T}(i,j)-2,0})$, so fixing some $R\geq 160M$ we see that
\[
d_{\cT{X}}(x,y)\leq d_{X}(\phi(x),\phi(y)) + 9R + \frac{1}{2}d_{\cT{X}}(x,y) + 2R + 16M.
\]
Hence,
\[
d_{\cT{X}}(x,y)\leq 2d_{X}(\phi(x),\phi(y)) + 22R + 32M.
\]

\section{Consequences of Theorem \ref{qTG}}\label{Conseq}
In this section we prove Corollaries $\ref{MCGtrees}$ and $\ref{RHdim}$. 

\subsection{Coned-off graphs and curve complexes}

We begin with one crucial result concerning embeddings of hyperbolic spaces into finite products of trees. Recall that a hyperbolic metric space $X$ is \textbf{visual} if for some (equivalently every) basepoint $x_0\in X$ there exists a $C>0$ and for each $x\in X$ a $(C,C)$-quasi--geodesic ray starting at $x_0$ and passing through $x$.

\begin{lem}\label{lem:visual} Let $X$ be a $\delta$-hyperbolic metric space with cobounded isometry group admitting a bi-infinite $(A,A)$ quasi-geodesic. Then $X$ is visual.
\end{lem}
\begin{proof} Fix some $x_0\in X$ and let $\gamma$ be a bi-infinite $(A,A)$ quasi-geodesic. There exists  a constant $D$ and for each $x\in X$ an isometry $\psi_x:X\to X$ such that $d_X(x,\psi(\gamma))\leq D$, so for each $x\in X$ there is a bi-infinite $(A,A+2D)$-quasi--geodesic $\gamma_x$ containing $x$.

Let $P_x$ be any path from $x_0$ to $\gamma_x$ with $\abs{P}\leq d_X(x_0,\gamma_x)+1$. Without loss we may assume $P_x\cap \gamma_x$ is a single point $x'$. Let $\gamma'_x$ be a $(A,A+2D)$ quasi-geodesic ray starting at $x'$ containing $x$. The concatenation $P_x \circ \gamma'_x$ is a $(C,C)$ quasi-geodesic ray starting at $x_0$ and passing through $x$, where $C$ depends only on $A,D$ and $\delta$.
\end{proof}

\begin{prop}\label{prop:asdimtotreeemb} Let $X$ be a visual hyperbolic space with asymptotic dimension at most $n$. Then $X$ admits a quasi-isometric embedding into a product of at most $n+1$ simplicial trees.
\end{prop}
\begin{proof} By \cite[Proposition $3.6$]{MS12} the capacity dimension of the boundary of $X$ is at most $n$, so the result follows from \cite[Theorem $1.1$]{Bu05b}.
\end{proof}

\begin{cor} Let $X$ be the curve complex of a compact surface, or a coned-off graph of a relatively hyperbolic group. Then $X$ admits a quasi--isometric embedding into a finite product of trees.
\end{cor}
\begin{proof} Curve complexes are hyperbolic and have finite asymptotic dimension \cite{MM99,BF08}, the mapping class group of the same surface acts coboundedly and the orbit of any pseudo-Anosov is a bi-infinite quasi--geodesic.

Similarly, the coned-off graph of a relatively hyperbolic group $G$ is hyperbolic, it has finite asymptotic dimension \cite[Theorem $5.1$]{Os05}, $G$ acts coboundedly on it, and the orbit of any element of infinite order which is not conjugate to an element in a peripheral subgroup is a bi-infinite quasi--geodesic.

In both cases the corollary follows from Lemma \ref{lem:visual} and Proposition \ref{prop:asdimtotreeemb}.
\end{proof}

To obtain explicit embeddings of curve complexes and coned-off graphs into $\ell^p$ spaces, we apply the methods of \cite{Hu11} to the collection of tight geodesics defined in \cite{Bo08} for curve complexes and in \cite{Bo08}[Corollary $3.5$] for coned-off graphs. The only properties we require from these papers are captured by the following definition (cf. \cite[Theorems $1.1$ and $1.2$]{Bo08}).

\begin{defn} Let $\Gamma$ be a graph. We say $\Gamma$ is a \textbf{Bowditch graph} if, for every pair of distinct vertices $a,b\in V\Gamma$ there is a set $\Gamma(a,b)$ containing a geodesic from $a$ to $b$ and such that the following conditions hold:
\begin{enumerate}
 \item For each $L$ there is a constant $K_0$ such that if $a,b\in V\Gamma$ and $c\in \Gamma(a,b)$ then $\Gamma(a,b)\cap B(c;L)$ has at most $K_0$ vertices.
 \item For each $L$ there exist constants $k_1$ and $K_1$ such that if $a,b\in V\Gamma$, $r\in\N$ and $c\in \Gamma(a,b)$ with $d(c,\set{a,b})\geq r+k_1$ then
 \[
   \bigcup_{d(a,x),d(b,y)\leq r} \Gamma(x,y) \cap B(c;L)
 \]
 contains at most $K_1$ vertices.
\end{enumerate}
\end{defn}

Notice that if $\Gamma$ is hyperbolic and has bounded geometry then it is automatically a Bowditch graph just by choosing $\Gamma(a,b)$ to be the set of all vertices lying on a geodesic from $a$ to $b$. More generally, we have the following, which appears in the comment after \cite[Corollary $3.5$]{Bo08}.

\begin{lem}\label{lem:unifineimpliesBowditch} Every uniformly fine hyperbolic graph is a Bowditch graph.
\end{lem}
Recall that a graph is uniformly fine if for every $n$ there is some $K(n)$ such that for any edge $e$ there are at most $K(n)$ simple loops of length $n$ containing $e$. Any fine graph admitting a group action with finitely many orbits of edges is uniformly fine, so the above result applies to coned-off graphs of relatively hyperbolic groups.
Lemma \ref{lem:unifineimpliesBowditch} is easily proved using the Morse property for geodesics in hyperbolic spaces.

\begin{thm}\label{a*cc} \h Let $X$ be a hyperbolic Bowditch graph and let $f:\N\to\N$ be a function satisfying the following
\[ f(n+1)-f(n)\leq f(n)-f(n-1)\tu{ for all }n\geq 1 \quad \tu{and} \quad\sum_{n\geq 1} \frac{1}{n}\left(\frac{f(n)}{n}\right)^p<\infty.
\]
Then there is an explicit embedding $\phi$ of $X$ into an $\ell^p$ space with
\[
f(d(x,y))\preceq \norm{\phi(x)-\phi(y)}_p \preceq d(x,y).
\]
In particular, $\alpha^*_p(X)=1$ for all $p\geq 1$.
\end{thm}
The second condition appears as property $(C_p)$ in \cite{Te11}.
\begin{proof} If one restricts attention to just the sets $\Gamma(a,b)$ considered in the above theorems of Bowditch then the result follows from carrying out the same procedure as in \cite[Section $2$]{Hu11}.
\end{proof}

\subsection{Embeddings of tree--graded spaces}

\begin{lem}\label{tgtoprodoftrees}
 Let $\mathcal{T}$ be a tree-graded space with pieces $\setcon{X_i}{i\in I}$. If, for each $i\in I$ there is a $(K,C)$ quasi-isometric embedding $q_i$ of $X_i$ into a product of $n$ trees $T_i^1\times\dots\times T_i^n$ then $\mathcal{T}$ quasi-isometrically embeds into a product of $n$ trees.
\end{lem}
\begin{proof} We assume the pieces are closed subsets of $\mathcal{T}$, if this is not the case we simply replace pieces by their closures.

Given a point $x\in \mathcal{T}$, let $\mathcal{T}_x$ be the ``transverse'' tree at $x$ constructed in \cite[Section $2$]{DS05}, namely it is the set of points which can be connected to $x$ by a topological arc whose intersection with any piece contains at most $1$ point. The subset $\mathcal{T}_x$ is an $\R$--tree \cite[Lemma $2.14(2)$]{DS05} and if $y\in \mathcal{T}_x$ then $\mathcal{T}_y=\mathcal{T}_x$.

We construct a tree-graded space $\mathcal{T}'$ with pieces $Y_i=T_i^1\times\dots\times T_i^n$ from a disjoint union of the $Y_i$ and the forest $\mathcal{F}=\bigsqcup \set{\mathcal{T}_x}$ by identifying $q_i(t)$ with $t$ whenever $t\in\mathcal{F}$.

Now $\mathcal{T}$ quasi-isometrically embeds into $\mathcal{T}'$ via the well-defined map $q(x)=q_i(x)$ whenever $x\in X_i$.

For each $j\in\set{1,\dots, n}$ define $T^j$ to be the space obtained from $\mathcal{T}'$ by projecting each piece $T_i^1\times\dots\times T_i^n$ onto $T_i^j$, label the projection $\mathcal{T}'\to T^j$ by $\pi_j$. $T^j$ is tree-graded with respect to pieces which are trees, so $T^j$ is a tree.

Now we claim that the map $\mathcal{T}'\to \prod_{j=1}^n T^j$  given by $x\mapsto (\pi_1(x),\dots,\pi_n(x))$ is a quasi-isometric embedding. For clarity, we denote the metric on $\mathcal{T}'$ by $d'$, the metric on each $T^j$ by $d_j$ and the metric on $\prod_{j=1}^n T^j$ by $d$. Since the $\pi_j$ are $1$-Lipschitz and, if we consider products with respect to an $\ell^1$ metric, we have
\[
 d'(x,y) \leq \sum_{j=1}^n d_j(\pi_j(x),\pi_j(y)) \leq n d(x,y)
\]
as required.
\end{proof}

\subsection{Proof of Corollary \ref{MCGtrees}}

Consider the surface $S=S_{g,n}$. If $3g+n-4\leq 0$ then $MCG(S)$ is virtually free and the results follow \cite{Be04}. We now assume $3g+n>4$.

Using the results of \cite{BBF10} together with Theorem $\ref{qTG}$ we obtain quasi-isometric embeddings of mapping class groups into finite products of tree-graded spaces, each of which have pieces uniformly quasi-isometric to curve complexes.
\[
 MCG(S)\to\prod_{i=1}^k\cC{\bY} \to \prod_{i=1}^k \cT{\bY}.
\]
There exist constants $K,C$ and $l$ such that each curve complex piece $\cC{Y}$ of $\cT{\bY}$ $(K,C)$--quasi-isometrically embeds into a product of $l$ regular trees of countable valence.

We therefore get a quasi-isometric embedding of $MCG(S)$ into a product of $kl$ trees using Lemma \ref{tgtoprodoftrees}.

Finally, we can construct an explicit embedding of a mapping class group into an $\ell^p$ space using Theorem \ref{a*cc} and \cite[Section $3$]{Hu11}.

\subsection{Proof of Corollary \ref{RHdim}}

The collection of maximal peripheral subgroups of a relatively hyperbolic group $G$, and near closest point projection maps with respect to a word metric on $G$, satisfies the axiomatic definition of \cite{BBF10}, so we obtain a space $\cC{\bH}$ together with an action of $G$ on this space. It follows from the projection description of relative hyperbolicity given by Sisto \cite[Theorem $0.1$]{Sis-projections} that the orbit map $G\to \cC{\bH}\times \hat{G}$ is a quasi-isometric embedding. Theorem $\ref{qTG}$ then implies that $G$ quasi-isometrically embeds into the product of a tree-graded space $\cT{\bH}$ with pieces uniformly quasi-isometric to subgroups $H_i$ with its coned-off graph $\hat{G}$.
\[
 G\to \cC{\bH}\times \hat{G} \to \cT{\bH}\times \hat{G}.
\]
Now $\hat{G}$ has finite Assouad-Nagata dimension, while $\cT{\bH}$ has this property if and only if each piece does; i.e. if the $H_i$ have finite Assouad-Nagata dimension.

To obtain an embedding into $\ell^1(\N)$ we use the fact that the coned-off graph quasi--isometrically embeds into a finite product of trees (and hence into $\ell^1(\N)$) and \cite[Section $3$]{Hu11}.

Finally, to obtain embeddings into $\ell^p$ spaces with the correct compression exponent, we use Theorem \ref{a*cc} and \cite[Section $3$]{Hu11}.

\end{document}